\documentclass{amsart}
\usepackage[margin=1.5in]{geometry}
\usepackage[colorlinks = true,linkcolor = blue,urlcolor  = blue,citecolor = blue,anchorcolor = blue]{hyperref}
\usepackage{amssymb,enumitem}
\usepackage{amsaddr}

\usepackage{xcolor,soul}

\newcommand{\R}{\mathbb R}

\newcommand{\N}{\mathbb N}

\newcommand{\E}{\mathbb E}

\renewcommand{\P}{\mathbb P}
\newcommand{\p}{\mathcal P}
\newcommand{\B}{\mathcal B}
\newcommand{\F}{\mathcal F}
\newcommand{\G}{\mathcal G}
\newcommand{\FF}{\mathbf F}
\newcommand{\s}{\mathcal S}
\newcommand{\M}{\mathcal M}
\newcommand{\A}{\mathcal A}
\newcommand{\1}{\mathbf{1}}

\newcommand{\bX}{\mathbf X}

\newcommand{\bigmid}{\,\middle\vert\,}
\newcommand{\reward}{\mathcal R}
\newcommand{\dd}{\mathrm{d}}

\newcommand{\esssup}{\operatorname*{ess\,sup}}
\newcommand{\essinf}{\operatorname*{ess\,inf}}

\theoremstyle{definition}

\newtheorem{theorem}{Theorem}[section]

\newtheorem{lemma}{Lemma}[section]
\newtheorem{definition}{Definition}[section]
\newtheorem{remark}{Remark}[section]
\newtheorem{proposition}{Proposition}[section]
\newtheorem{example}{Example}[section]

\numberwithin{equation}{section}

\begin{document}

\title[Conditional independence and the Gittins index]{The Gittins index is optimal for dynamic allocation with conditionally independent filtrations}

\author[C.~Wang]{Christopher Wang}
\thanks{This research was supported in part by the NSF GRFP under Grant No.~DGE--2139899.
The author would like to express immense thanks to Dr.~Ioannis Karatzas for proposing this topic, in addition to his suggestions, insight, and continued support; as well as the anonymous referees whose comments have greatly improved the contents and structure of this paper.}
\address{\small Department of Mathematics, Cornell University, Ithaca, NY 14853, USA}
\email{cyw33@cornell.edu}

\begin{abstract}
The dynamic allocation problem, also known as the `multi-armed bandit' problem, simulates a situation in which an agent is faced with a tradeoff between actions that yield an immediate reward and actions whose benefits can only be perceived in the future. In this paper, we show that the non-Markovian, discrete-time problem can be solved by following a Gittins index strategy, without the assumption that the rewards processes are independent. Instead, we require the underlying multi-parameter filtration to satisfy a conditional independence property. We provide three representations of the maximal attainable value under an optimal strategy. Furthermore, we discuss the relationship between index-type strategies and the `synchronization' paradigm from operations research.
\end{abstract}



\maketitle
\thispagestyle{empty}

\section{Introduction}\label{s:intro}

The dynamic allocation problem, also known as the `multi-armed bandit' problem, simulates a situation in which an agent is faced with a tradeoff between actions that yield an immediate reward and actions whose benefits can only be perceived in the future. In a classical example, an investigator is faced with $d$ `projects' which simultaneously compete for attention. At any given discrete moment in time, only one of these projects may be engaged while the others remain `frozen.' Upon engaging a project, the investigator receives a random reward, which becomes discounted over time. Such a scenario encapsulates what is known in machine learning as the conflict between `exploitation' and `exploration.' On the one hand, the investigator may choose to consistently engage a small subset of well-understood projects whose rewards are acceptably high; on the other hand, the investigator may find it worthwhile to spend time `learning' about each of the projects, perhaps sacrificing some greater immediate rewards in doing so, to devise a stronger strategy for the future. A solution to the dynamic allocation problem consists of finding the sequence of projects that should be engaged which maximizes the overall rewards reaped by the investigator, as well as the value of the rewards acquired through such an optimal strategy.

To be precise, for each moment in time $t\in\N_0$, let $T_i(t)$ denote the total number of occasions one has engaged the $i$th project by time $t$, where $i$ ranges from $1,\dots,d$. We require that $\sum_{i=1}^d T_i(t) = t$, and that if $T_i(t+1) = T_i(t)+1$, then $T_j(t)=T_j(t+1)$ for each $j\neq i$, encapsulating the constraint that one may only engage a single project at a time, during which the remaining projects are frozen. By engaging project $i$ at time $t$, one receives a reward $h_i(T_i(t))$, geometrically discounted in time by a factor of $\beta\in(0,1)$. Finally, one may abandon the entire endeavor in return for a `retirement reward' $M>0$, also discounted with time. The dynamic allocation problem is to determine an optimal `allocation strategy' $\tilde T^*(t) = (T_1^*(t),\dots,T_d^*(t))$, $t\in\N_0$, along with an optimal stopping time $\eta^*$, which attains the maximum total expected reward
$$\Phi := \sup_{\tilde T,\eta}\,\E\left[\sum_{i=1}^d \sum_{t=0}^{\eta-1} \beta^t h_i(T_i(t+1)) \Big(T_i(t+1)-T_i(t)\Big) + M\beta^\eta\right],$$
and to compute the value of $\Phi$. We assume the rewards sequences $H_i:=\{h_i(t)\}_{t\in\N_0}$ are positive, uniformly bounded, and adapted to filtrations $\FF_i:=\{\F_i(t)\}_{t\in\N_0}$ representing their histories; we do not require them necessarily to be Markov processes. We do however require that decisions about the allocation strategy and the retirement time be made only on the basis of information that has already been accumulated in the course of the investigation. In other words, we require $\tilde T$ and $\eta$ be adapted to a so-called `multi-parameter filtration' $\FF:=\{\F(\tilde s)\}_{\tilde s\in\N_0^d}$, related to the project filtrations $\FF_i$, which represents the collective history of all $d$ projects accumulated up to moments $\tilde s$ in `multi-parameter time.'

The multi-armed bandit problem enjoys a long history of research, beginning as early as World War II, during which ``efforts to solve it so sapped the energies and minds of Allied analysts that the suggestion was made that the problem be dropped over Germany, as the ultimate instrument of intellectual sabotage,'' in the words of Peter Whittle \cite{Discussion1979}. Building on early work in experimental design from the 1920s, Herbert Robbins outlined the first modern formulation of the problem in 1952 \cite{Robbins1952}, yet the first major breakthrough did not take place until 1979 with Gittins' brilliant development of the so-called `Gittins index' and the index-type strategy \cite{Gittins1979}. Gittins' key insight was to introduce an additional parameter to the problem by allowing the investigator to abandon, at any point in time, the entire collection of projects in exchange for a fixed `retirement reward.' The effect of this new parameter, which {\em prima facie} makes the problem even more difficult, is that it allows one to recognize each of the $d$ projects as its own {\em optimal stopping problem}, which can be solved independently and `stitched together' to form an optimal strategy for the dynamic allocation problem. Additional insights in both discrete and continuous time for Markov rewards processes followed in the 1980s, e.g. \cite{Whittle1980, Karatzas1984, Tsitsiklis1986}. 

Theory for the non-Markovian setting, in which the random rewards additionally depend on the history of the project being engaged, was first developed in \cite{Varaiya1985} and afterward split largely into two camps: the `excursion-theoretic' approach of \cite{Kaspi1994, Kaspi1998}, and the `martingale-theoretic' approach of \cite{Mandelbaum1986, Mandelbaum1987, ElKaroui1993, ElKaroui1994, ElKaroui1997}. With the exception of \cite{ElKaroui1997}, all formulations of dynamic allocation assumed that the projects evolved {\em independently} in time, i.e., that the choice to engage a given project has no effect on the rewards yielded by other projects. In 1997, \cite{ElKaroui1997} showed in the continuous-time case that there exists an optimal index-type strategy even when the projects were not assumed to be independent. They instead assume only the celebrated conditional independence condition, frequently referred to as the \eqref{itm:F4} condition of \cite{Cairoli1975}, which was originally developed for multi-parameter martingale theory. Heuristically, the \eqref{itm:F4} condition relaxes the independence assumption by allowing the different projects to be dependent upon a `shared history,' whereas they become independent conditioned on one's knowledge of that shared history.

In this work, we formulate and extend the results of \cite{ElKaroui1997} to the discrete-time setting, and we discuss additional properties of index-type strategies and of representations of the dynamic allocation problem under the \eqref{itm:F4} condition. One recurring finding is that the random rewards sequence can be replaced by a certain {\em almost surely decreasing} sequence of rewards, without affecting the expected maximum achievable rewards or the optimal strategy. This fact was previously known  in the one-dimensional case (Equation \eqref{eq:ost-decr-rep-2}) and in the case where the projects are independent \cite{ElKaroui1993a, ElKaroui1993}; the main result of this work, Theorem \ref{thm:main}, shows that this remains true in the most general scenario, where the filtrations of the multiple projects satisfy only the \eqref{itm:F4} condition and are not necessarily independent. The proof of Theorem \ref{thm:main} relies on a class of strategies satisfying the so-called `synchronization identity,' which strictly include those of index type. 

By working in the discrete-time setting, we make several new contributions to the theory of dynamic allocation:
\begin{enumerate}
\item We establish the equivalence of three different representations of the value of the dynamic allocation problem (Theorems \ref{thm:main} and \ref{thm:dap-indep}); these equivalences were previously confirmed only when assuming independence of filtrations, and unknown to hold in continuous time under the weaker (F4) condition.
\item Compared to \cite{ElKaroui1997}, which constructs only one strategy that is optimal for the dynamic allocation problem in the (F4) setting, we show that their strategy is in fact of index type (Proposition \ref{prop:y-strategy}), and that all index-type strategies are optimal.
\item We fully characterize, in discrete time, the relationship between the class of strategies satisfying the synchronization identity and the class of index-type strategies, by means of excursion theory (Proposition \ref{prop:excursion}).
\end{enumerate}
Additionally, we have aimed to provide detailed and clear proofs throughout, and with the exception of some results in the theory of optimal stopping and multi-parameter martingale theory, this paper is essentially self-contained. Our opinion is that the discrete-time setting permits a clarified presentation of the results of \cite{ElKaroui1993,ElKaroui1997}, upon which this work is largely built. In particular, we are able to give relatively elementary proofs of time-change results without relying on the machinery of stochastic calculus.

Finally, we remark that the dynamic allocation problem in general is well-known to have many important applications in a wide range of fields, such as experimental design for clinical trials \cite{Villar2015}, path routing in communication networks \cite{Awerbuch2008}, and algorithms for machine learning \cite{Radlinski2008}. For potential applications of the non-Markovian setup under a conditional independence assumption on the filtrations, \cite[Sect.~5.1]{Scully2025} provides an example of a cost-per-sample Bayesian optimization problem with an index strategy that incorporates correlation data from prior samples. Other recent applications of non-Markovian dynamic allocation include reinforcement learning settings in which tasks must be performed in a specific sequence in order to reap rewards \cite{Gaon2020}, as well as the optimal scheduling of healthcare interventions in a manner that incorporates empirically observed temporal dependencies in patients' behavior \cite{Danassis2023}.

We begin in Section \ref{s:osp} with a brief overview of the one-dimensional setting, which is simply an optimal stopping problem, and we introduce Gittins index sequences in Section \ref{s:gittins}, following \cite{ElKaroui1993}. In Section \ref{s:dap}, we introduce the dynamic allocation problem in full generality and determine its value and optimal strategy in the setting where the reward sequences are almost surely decreasing. In Section \ref{s:strategies}, we detail the implications of the synchronization identity  and delineate the relation between the synchronization identity and strategies of index type by means of excursion theory. Multi-parameter martingale theory is reviewed in Section \ref{s:multiparam}. Finally, results from each of the previous sections are used in Section \ref{s:f4} to state and prove our main theorem, Theorem \ref{thm:main}, and we conclude in Section \ref{s:whittle} with Whittle's representation of the value of the dynamic allocation problem.

%
%

\section{The optimal stopping problem}\label{s:osp}

We consider first the problem of optimal stopping in discrete time, equivalent to dynamic allocation with retirement reward when $d=1$. Proofs for results in this section can be found in \cite{Neveu1975}, unless otherwise stated.

Let us work on the complete probability space $(\Omega,\F,\P)$ equipped with filtration $\FF = \{\F(t)\}_{t\in\N_0}$. Let $H = \{h(t)\}_{t\in\N}$ be the positive `rewards sequence,' which
\begin{enumerate}
\item is predictable with respect to $\FF$, meaning that $h(t)$ is $\F(t-1)$-measurable, for every $t\in\N$; and
\item satisfies the integrability condition $\E[\sum_{t=0}^\infty\beta^t h(t+1)]<\infty,$
\end{enumerate}
where $\beta\in(0,1)$ is the `discount factor.' Let $\s(\theta)$ denote the collection of
$\FF$-stopping times taking value in $\{\theta, \theta+1, \dots\}\cup\{+\infty\}$ for some $\theta\in\N_0$.

\begin{definition}
For `exit reward' $m\in[0,\infty)$ and time $t\in\N_0$, the {\bf optimal stopping problem} is to compute the value
\begin{equation}\label{eq:osp}
V(t;m) := \esssup_{\tau\in\s(t)}\,\E\left[\sum_{u=t}^{\tau-1} \beta^{u-t} h(u+1) 
	+ m\beta^{\tau-t} \bigmid \F(t)\right].
\end{equation}
Here, $V(t;m)$ is the expected discounted reward indexed by the exit reward $m$, given that the agent begins at time $t$ and possesses the information available at that time. If a stopping time $\tau\in\s(t)$ attains this supremum, we call it an {\bf optimal stopping time}. Notice that taking $\tau\equiv t$ immediately yields
\begin{equation}\label{eq:V-lower}
V(t;m)\ge m, \qquad \text{a.s.}
\end{equation}
for all $t\in\N_0$, $m\ge0$.
\end{definition}

Recall the following processes, standard in the theory of optimal stopping. We denote by
\begin{equation}\label{eq:Y-def}
Y(t;m) := \sum_{u=0}^{t-1} \beta^u h(u+1) + m\beta^t, \qquad t\in\N
\end{equation}
the total reward received upon stopping at time $t$ given exit reward $m\ge0$. The `Snell envelope' of $Y(\,\cdot\,;m)$ is given by
\begin{equation}\label{eq:snell}
Z(t;m) := \esssup_{\tau\in\s(t)}\,\E[Y(\tau;m)\mid\F(t)], \qquad t\in\N_0
\end{equation}
and satisfies
\begin{equation}\label{eq:Z-1}
Z(t;m) = \sum_{u=0}^{t-1} \beta^u h(u+1) + \beta^t V(t;m).
\end{equation}
The random sequence $\{Z(t;m)\}_{t\in\N_0}$ is the smallest $\FF$-supermartingale that dominates the rewards sequence $\{Y(t;m)\}_{t\in\N}$, hence the terminology of `envelope.' It happens that the relationship between $Z$ and $Y$ determines the optimal stopping time \cite{Snell1952}.
	
\begin{theorem}\label{thm:ost-soln}
For the processes $Y,Z$ as defined in \eqref{eq:Y-def}, \eqref{eq:snell}, the stopping time
\begin{equation}\label{eq:sigma}
\sigma(t;m) := \inf\{\theta\ge t\colon Y(\theta;m) = Z(\theta;m)\}
\end{equation}
is contained in $\s(t)$ and is optimal for problem \eqref{eq:osp}. In other words,
\begin{equation}\label{eq:sigma-optimal}
V(t;m) = \E\left[\sum_{u=t}^{\sigma(t;m)-1} \beta^{u-t}h(u+1) + m\beta^{\sigma(t;m)-t} 
	\bigmid \F(t)\right].
\end{equation}
\end{theorem}
From \eqref{eq:Z-1} and \eqref{eq:sigma}, we get the natural relationship
\begin{equation}\label{eq:sigma-V}
\sigma(t;m) = \inf\{\theta\ge t \colon V(\theta;m) = m\},
\end{equation}
which says that it is optimal to stop as soon as the future expected discounted reward, conditioned on the currently available information, is exactly equal to the reward the agent would receive from exiting immediately.
	
From the theory of optimal stopping, we also have the following results.
	
\begin{proposition}\label{prop:Z-mtg}
The stopped process $\{Z(\theta\wedge\sigma(t;m);m), \F(\theta)\}_{\theta=t}^\infty$ is a martingale for every $t\in\N_0$ and $m\ge0$.
\end{proposition}

\begin{proposition}[Dynamic Programming Equation]
For every $t\in\N_0$ and $m\ge0$,
\begin{equation}\label{eq:dp}
V(t;m) = \max\Big(m, h(t+1) + \beta\,\E[V(t+1;m)\,\big|\,\F(t)]\Big)
\end{equation}
holds almost surely.
\end{proposition}
	
The dynamic programming equation states that the optimal future expected discounted reward can be achieved either by exiting immediately and collecting the exit reward $m$, or by continuing for one additional unit of time and evaluating the optimal future conditional expected reward at that point in time.
	
Finally, we relate two lemmas from \cite{ElKaroui1993} which we require later, with self-contained proofs in the \hyperref[appendix]{Appendix}.
		
\begin{lemma}\label{lem:sigma-props}
For any $t\in\N_0$, the map $m\mapsto\sigma(t;m)$ is almost surely decreasing and right-continuous.
\end{lemma}
	
\begin{lemma}\label{lem:V-rdiff}
For any $t\in\N_0$, the map $m\mapsto V(t;m)$ is convex, increasing, and has right-derivative
\begin{equation}\label{eq:V-rdiff}
\frac{\partial^+}{\partial m}V(t;m) := \lim_{\delta\downarrow0}\frac{V(t;m+\delta) - V(t;m)}
	{\delta} = \E\left[\beta^{\sigma(t;m)-t}\mid\F(t)\right], \qquad \text{a.s.}
\end{equation}
\end{lemma}
	
	%
	%
	
\section{Gittins index sequences}\label{s:gittins}
	
In the dynamic allocation problem with retirement reward, the investigator always has the option of disregarding all projects but one; by focusing on a single project and excluding all others, the investigator can treat each project as its own optimal stopping problem and construct an allocation strategy based on the project's `Gittins index.' We define the Gittins index for each such optimal stopping problem in isolation from the others; results in this section generally follow \cite{ElKaroui1993}.
	
\begin{definition}
Let $\M(t)$ be the family of positive $\F(t)$-measurable random variables for $t\in\N_0$. Define the {\bf Gittins index} at time $t$ by
\begin{equation}\label{eq:gittins}
M(t) := \esssup\{X\in\M(t)\colon V(t;X) \ge X\} = \essinf\{X\in\M(t)\colon V(t;X) = X\},
\end{equation}		
which represents the minimum reward the investigator must be offered to exit immediately, at the current time $t$, i.e., `the equitable surrender value.' We also define the {\bf minimum Gittins index to date}, or the `lower envelope' of the Gittins index to date, for given $t,\theta\in\N_0$ with $t\le\theta$:
\begin{equation}\label{eq:gittins-envelope}
\underline M(t,\theta) := \min_{t\le u\le\theta} M(u), \qquad \qquad
\underline M(\theta) := \underline M(0,\theta).
\end{equation}
\end{definition}
	
\begin{remark}
The definition \eqref{eq:gittins} of the Gittins index is due to \cite{Whittle1980}, whereas Gittins' original formulation is via the `forward induction' characterization \cite[Eq.~(4.1)]{Gittins1979}:
\begin{equation}\label{eq:gittins-forwards}
(1-\beta)M(t) = \esssup_{\tau\in\s(t+1)} 
	\frac{\E\left[\sum_{u=t}^{\tau-1}\beta^uh(u+1)\bigmid\F(t)\right]}
	{\E\left[\sum_{u=t}^{\tau-1}\beta^u\bigmid\F(t)\right]}, \qquad\text{a.s.}
\end{equation}
Proofs for the equivalence of the two formulations are given in \cite{Mandelbaum1986,ElKaroui1993a}; we present a new method in the \hyperref[appendix]{Appendix}.
\end{remark}
	
Clearly the mapping $\theta\mapsto\underline M(t,\theta)$ is decreasing. The next result demonstrates that the mappings $\theta\mapsto \underline M(t,\theta)$ and $m\mapsto\sigma(t;m)$ are, in a sense, right-inverses of each other. As a consequence, understanding the behavior of the Gittins index over time will reveal the behavior of the optimal stopping rule as a function of the exit reward. 
	
\begin{lemma}\label{lem:sigma-equivs}
For $m\ge0$ and $t,\theta\in\N_0$ with $t\le\theta$, we have almost surely that
\begin{equation}\label{eq:sigma-equivs}
\sigma(t;m) > \theta\ \iff\ V(s;m) > m,\ \forall\,s=t,\dots,\theta\ \iff\ 
	m < \underline M(t,\theta).
\end{equation}
\end{lemma}
\begin{proof}
The first equivalence follows directly from \eqref{eq:V-lower} and \eqref{eq:sigma-V}. Now we provethe second equivalence.
		
Assume that $V(s;m)>m$ holds almost surely for $s=t,\dots,\theta$, and suppose for some $s$ that $M(s)=m$ holds on a set $A\subset\F$ of positive measure. Recall from \eqref{eq:osp} that the map $V(s;\cdot)$ is continuous, hence the set $\{\tilde m\ge0\colon V(s;\tilde m)(\omega)>\tilde m\}$ is open, and thus does not contain its supremum, for almost every $\omega\in\Omega$. It follows from from \eqref{eq:gittins} that for almost every $\omega\in A$, we have
$$M(s)(\omega) = m = \sup\{\tilde m\ge0\colon V(s;\tilde m)(\omega)>\tilde m\} \in
	\{\tilde m\ge0\colon V(s;\tilde m)(\omega)\le\tilde m\}.$$
Thus $V(s;m)(\omega)\le m$ holds on a set of positive measure, contradicting our initial assumption, and so $M(s)>m$ holds almost surely, for every $s=t,\dots,\theta$. We now obtain $\underline M(t,\theta)>m$, as desired, from definition \eqref{eq:gittins-envelope}.
		
Conversely, $\underline M(t,\theta) > m$ gives $M(s)>m$ for every $s=t,\dots,\theta$. We know from \eqref{eq:V-lower} that $V(s;M(s))\ge M(s)$ holds almost surely. Now recall from \eqref{phi decrease} that the map $\tilde m\mapsto V(s;\tilde m)-\tilde m$, is strictly decreasing, hence 
$$V(s;m)-m>V(s;M(s))-M(s)\ge0, \qquad \text{a.s.}$$
holds for each $s=t,\dots,\theta$.
\end{proof}
	
Since $\sigma(t;m)$ is an integer, $\sigma(t;m)=\theta$ holds if and only if $\sigma(t;m)>\theta-1$ and $\sigma(t;m)\le\theta$, which by Lemma \ref{lem:sigma-equivs} is equivalent to $\underline M(t,\theta-1)>m $ and $\underline M(t,\theta)\le m$. Since $\theta\mapsto\underline M(t,\theta)$ is decreasing, we obtain
\begin{equation}\label{eq:sigma-M}
\sigma(t;m) = \inf\{\theta\ge t\colon \underline M(t,\theta)\le m,\ \text{a.s.}\}
	= \inf\{\theta\ge t\colon  M(\theta)\le m,\ \text{a.s.}\}.
\end{equation}
Analogously, the right-continuity of $m\mapsto\sigma(t;m)$ yields
\begin{equation}\label{eq:M-sigma}
\underline M(t,\theta) = \inf\{m\ge0\colon \sigma(t;m)\le\theta,\ \text{a.s.}\} 
	= \sup\{m\ge0\colon \sigma(t;m)>\theta,\ \text{a.s.}\}.
\end{equation}
	
\begin{proposition}
For $t,\theta\in\N_0$ with $t\le\theta$, we have
\begin{equation}\label{eq:sigma-jumps}
\sigma(t;\underline M(t,\theta))\le \theta < \sigma(t;\underline M(t,\theta)-),\qquad\text{a.s.}
\end{equation}
\end{proposition}
\begin{proof}
The inequalities that follow hold almost surely. By \eqref{eq:M-sigma}, there exists a sequence $m_n\downarrow\underline M(t,\theta)$ such that $\sigma(t;m_n)\le\theta$ holds for each $n\in\N$, so the right-continuity of $m\mapsto\sigma(t;m)$ yields $\sigma(t;\underline M(t,\theta))\le\theta$.
		
On the other hand, there exists a sequence $m'_n\uparrow\underline M(t,\theta)$, also by \eqref{eq:M-sigma}, such that $\sigma(t;m'_n)>\theta$ holds for each $n\in\N$. Since $\sigma(t,\cdot)$ takes integer values we obtain $\sigma(t;m'_n)\ge\theta+1$, which implies $\sigma(t;\underline M(t,\theta)-)\ge\theta+1>\theta$.
\end{proof}
	
\begin{proposition}\label{prop:M-constant}
If $\sigma(t;\,\cdot\,)$ has a jump at $m\ge0$ for fixed $t\in\N_0$ then the mapping $\theta\mapsto\underline M(t,\theta)$ is almost surely constant on $\{\sigma(t;m),\dots,\sigma(t;m-)-1\}$, and
\begin{equation}\label{eq:M-constant}
m = \underline M(t,\sigma(t;m)) = \underline M(t,\sigma(t;m)+1)
	= \dots = \underline M(t,\sigma(t;m-)-1),\qquad\text{a.s.}
\end{equation}
\end{proposition}
\begin{proof}
Fix an integer $\sigma(m)\le\theta < \sigma(m-)$. By \eqref{eq:sigma-equivs}, $\sigma(t,m)\le\theta$ implies that $m\ge\underline M(t,\theta)$ holds almost surely. On the other hand, the decrease of $\sigma(t;\cdot)$ implies that $\sigma(t;m-\delta)\ge\sigma(t;m-)>\theta$ holds almost surely for every $0<\delta\le m$. Thus, $m-\delta<\underline M(t,\theta)$ holds almost surely, again by \eqref{eq:sigma-equivs}, and taking $\delta\downarrow0$ yields $m\le\underline M(t,\theta)$, almost surely. 
\end{proof}
	
From the relation between $\sigma(t;\,\cdot\,)$ and $\underline M(t,\,\cdot\,)$, we deduce several notable properties of the optimal stopping problem.
	
\begin{lemma}\label{lem:V-rdiff-props}
For any $t\in\N_0$, the mapping $m\mapsto\frac{\partial^+}{\partial m}V(t;m)$ is almost surely right-continuous and increasing. Moreover, it takes the value zero for $m=0$ and the value one for $m\ge M(t)$, almost surely.
\end{lemma}
\begin{proof}
Right-continuity and increase follow immediately from Lemma \ref{lem:sigma-props} and \eqref{eq:V-rdiff}. By \eqref{eq:gittins} and \eqref{eq:sigma-equivs}, we have $M(\theta)>0$ and thus $\sigma(t;0)>\theta$ for all $\theta\ge t$. Therefore 
\begin{equation}\label{eq:sigma-0}
\sigma(t;0)=+\infty, \qquad\text{a.s.},
\end{equation}
which in conjunction with \eqref{eq:V-rdiff} implies $\frac{\partial^+}{\partial m}V(t;0)=0$. On the other hand, if $M(t)\le m$ then \eqref{eq:sigma-M} yields $\sigma(t;m)=t$, and so $\frac{\partial^+}{\partial m}V(t;m)=1$.
\end{proof}
	
We also obtain a representation of the value of the optimal stopping problem in terms of the Gittins indices. Note the significance of \eqref{eq:ost-decr-rep-2} below: when computing future expected discounted rewards, conditioned on current information, the rewards sequence $\{h(\theta+1)\}_{\theta=t}^\infty$ can be replaced by the \emph{decreasing} sequence $\{(1-\beta)\underline M(t,\theta)\}_{\theta=t}^\infty$, without consequence. 
	
\begin{proposition}\label{prop:ost-decr-rep}
For any $t\in\N_0$ and $m\ge0$, we have the almost sure representation
\begin{multline}\label{eq:ost-decr-rep-1}
V(t;m) = (1-\beta)\,\E\left[\sum_{\theta=t}^\infty \beta^{\theta-t}
	\Big(m\vee \underline M(t,\theta)\Big) \bigmid \F(t)\right] \\
= m + \E\left[\int_m^\infty \left(1-\beta^{\sigma(t;\lambda)-t}\right)\,\dd\lambda\bigmid \F(t)\right]
\end{multline}
of the value in the optimal stopping problem \eqref{eq:osp}, in terms of the minimum Gittins index to date in \eqref{eq:gittins-envelope}.
\end{proposition}
\begin{proof}
First notice that \eqref{eq:sigma-optimal} and \eqref{eq:sigma-0} yield
\begin{align}\label{eq:ost-1}
V(t;m) - V(t;0) = \E\left[m\beta^{\sigma(t;m)-t} 
	- \sum_{\theta=\sigma(t;m)}^\infty\beta^{\theta-t}h(\theta+1) \bigmid \F(t)\right].
\end{align}
Now, for every $\tilde m\le m$, we write $\beta^{\sigma(t;\tilde m)-t}$ as a geometric sum and apply the equivalence \eqref{eq:sigma-equivs} to obtain
\begin{align*}
\beta^{\sigma(t;\tilde m)-t}
&= (1-\beta)\sum_{\theta=\sigma(t;m)}^\infty\beta^{\theta-t}
	\1\{\underline M(t,\theta)\le \tilde m\}.
\end{align*}
Substituting this value into \eqref{eq:V-rdiff} and integrating yields
\begin{align*}
V(t;m) - V(t;0)
&= \int_0^m\frac{\partial^+}{\partial\tilde m}V(t;\tilde m)\,\dd\tilde m \\
&= (1-\beta)\int_0^m \E\left[\sum_{\theta=\sigma(t;m)}^\infty\beta^{\theta-t}
	\1\{\underline M(t,\theta)\le \tilde m\} \bigmid \F(t)\right] \dd\tilde m,
	\qquad\text{a.s.}
\end{align*}
Note that the fundamental theorem of calculus can be used to establish the first equality above, despite the fact that the right-derivative of $V(t;\cdot)$ is used; this is due to convexity, established in Lemma \ref{lem:V-rdiff} (e.g. \cite[Prob.~3.6.21]{Karatzas1988}). Since the integrand on the right-hand side is non-negative, Tonelli's theorem implies
\begin{equation}\label{eq:ost-2}
V(t;m)-V(t;0) = (1-\beta)\,\E\left[\sum_{\theta=\sigma(t;m)}^\infty
	\beta^{\theta-t}\Big(m-\underline M(t,\theta)\Big) \bigmid \F(t)\right], \qquad \text{a.s.},
\end{equation}
where we use the fact that $m-\underline M(t,\theta)\ge0$ when $\theta\ge\sigma(t;m)$, from \eqref{eq:sigma-equivs}.
		
Finally, we note
\begin{align}\label{eq:ost-3}
(1-\beta)&\,\E\left[\sum_{\theta=\sigma(t;m)}^\infty
	\beta^{\theta-t} \Big(m-\underline M(t,\theta)\Big) \bigmid \F(t)\right] \nonumber \\
&= m\,\E[\beta^{\sigma(t;m)-t} \mid \F(t)]
	+ (1-\beta)\,\E\left[\sum_{\theta=\sigma(t;m)}^\infty
	\beta^{\theta-t}\underline M(t,\theta) \bigmid \F(t)\right].
\end{align}
Recalling \eqref{eq:ost-2}, we deduce that the right-hand sides of \eqref{eq:ost-1} and \eqref{eq:ost-3} are almost surely equal, hence
\begin{equation}\begin{split}\label{eq:ost-4}
\E\left[\sum_{\theta=\sigma(t;m)}^\infty\beta^{\theta-t}h(\theta+1) \bigmid \F(t)\right]
	= (1-\beta)\, \E\left[\sum_{\theta=\sigma(t;m)}^\infty
	\beta^{\theta-t}\underline M(t,\theta) \bigmid \F(t)\right],\quad\text{a.s.}
\end{split}\end{equation}
Once more, \eqref{eq:sigma-equivs} implies that $\sigma(t;m)\to t$ as $m\to\infty$, so it follows from \eqref{eq:sigma-optimal}, \eqref{eq:sigma-0}, and \eqref{eq:ost-4} that
\begin{equation}\label{eq:ost-decr-rep-2}
V(t;0) = \E\left[\sum_{\theta=t}^\infty\beta^{\theta-t}h(\theta+1) \bigmid \F(t)\right]
	= (1-\beta)\, \E\left[\sum_{\theta=t}^\infty
	\beta^{\theta-t}\underline M(t,\theta) \bigmid \F(t)\right],\quad \text{a.s.}
\end{equation}
By substituting the expression on the right-hand side into \eqref{eq:ost-2}, we obtain the desired representation of $V(t;m)$.

The second equality follows from \eqref{eq:sigma-equivs}, namely
\begin{align*}
(1-\beta)&\E\left[\sum_{\theta=t}^\infty 
	\beta^{\theta-t}(m\vee\underline M(t,\theta))\bigmid\F(t)\right] \\
&\qquad= (1-\beta)\E\left[\sum_{\theta=t}^\infty \beta^{\theta-t}\left(m + \int_m^\infty\1\{\sigma(t;\lambda)>\theta\}\,\dd\lambda\right)\bigmid\F(t)\right] \\
&\qquad= m + \E\left[\int_m^\infty\left(1-\beta^{\sigma(t;\lambda)-t}\right)\,\dd\lambda\bigmid\F(t)\right]
\end{align*}
as desired.
\end{proof}
	
The final result in this section now compares the expected rewards, evolving in time, when receiving the actual rewards $h(\cdot)$ versus the decreasing rewards $(1-\beta)\underline M(t,\cdot)$.
	
\begin{proposition}\label{prop:U-mtg}
For every $t\in\N_0$ and $m\ge0$, the sequence
\begin{equation*}
\left\{U(\theta) \vphantom{\sum_{u=t}^{\theta-1}} := \beta^\theta
	\Big[V(\theta;\underline M(t,\theta)) - \underline M(t,\theta)\Big] 
	+ \sum_{u=t}^{\theta-1} \beta^u\Big(h(u+1)-(1-\beta)\underline M(t,u)\Big)
	\right\}_{\theta=t}^\infty 
\end{equation*}
is an $\{\F(\theta)\}_{\theta=t}^\infty$-martingale.
\end{proposition}
\begin{proof}
The equalities below hold almost surely.
		
Fix $\theta\ge t$. We recall that $(1-\beta)^{-1}=\sum_{u=0}^\infty\beta^u$, hence \eqref{eq:ost-decr-rep-1} implies
\begin{equation}\label{eq:U-1}
\beta^\theta\Big[V(\theta;m)-m\Big] = (1-\beta)\,\E\left[\sum_{u=\theta}^\infty 
	\beta^u\Big(\underline M(\theta,u)-m\Big)^+\bigmid\F(\theta)\right].
\end{equation}
This holds for every $m\ge0$, so replacing $m$ by $\underline M(t,\theta)$, an $\F(\theta)$-measurable random variable, yields
\begin{equation}\label{eq:U-2}
\beta^\theta\Big[V(\theta;\underline M(t,\theta))-\underline M(t,\theta)\Big]
	= (1-\beta)\,\E\left[\sum_{u=\theta}^\infty 
	\beta^u\Big(\underline M(\theta,u)-\underline M(t,\theta)\Big)^+\bigmid\F(\theta)\right].
	\end{equation}
Now fix $m\ge0$. On the event $\{\theta<\sigma(t;m)\}=\{m<\underline M(t,\theta)\}$, we easily check that
\begin{align*}
\Big(\underline M(\theta,u)-m\Big)^+
	- \Big(\underline M(\theta,u)-\underline M(t,\theta)\Big)^+
&= \Big(\underline M(t,u)-m\Big)^+
\end{align*}
for every $u\ge\theta$. Thus, by subtracting \eqref{eq:U-2} from \eqref{eq:U-1}, we obtain
\begin{equation}\label{eq:U-3}
\beta^\theta D(\theta) = (1-\beta)\,\E\left[\sum_{u=\theta}^\infty 
	\beta^u\Big(\underline M(t,u)-m\Big)^+\bigmid\F(\theta)\right]
\end{equation}
on $\{\theta<\sigma(t;m)\}$, where
\begin{equation}\label{eq:U-4}
D(\theta) := \Big(V(\theta;m)-m\Big) - \Big(V(\theta;\underline M(t,\theta)) 
	- \underline M(t,\theta)\Big),\qquad\theta\ge t.
\end{equation}
		
For ease of notation, let $\theta_\sigma(m) := \theta\wedge\sigma(t;m)$. We claim that the sequence
\begin{equation}\label{eq:U-5}
\left\{\beta^{\theta_\sigma(m)} D(\theta_\sigma(m))
	+ (1-\beta)\sum_{u=t}^{\theta_\sigma(m)-1}\beta^u(\underline M(t,u)-m),
	\right\}_{\theta=t}^\infty
\end{equation}
is an $\{\F(\theta)\}_{\theta=t}^\infty$-martingale. Indeed, fix $\theta,\theta'\in\N_0$ such that $t\le\theta\le\theta'$. Then
\begin{align*}
\E\left[ \vphantom{\sum_{u=t}^{\theta'_\sigma(m)-1}} \beta^{\theta'_\sigma(m)} \right.
	&\left. D(\theta'_\sigma(m)) + (1-\beta)\sum_{u=t}^{\theta'_\sigma(m)-1}
	\beta^u\Big(\underline M(t,u)-m\Big)\bigmid \F(\theta)\right] \\
&= (1-\beta)\, \E\left[\E\left[\sum_{u=t}^\infty\beta^u\Big(\underline M(t,u)-m\Big)^+
	\bigmid\F(\theta'_\sigma(m))\right]\bigmid\F(\theta)\right] \\
&= (1-\beta)\,
	\E\left[\sum_{u=\theta_\sigma(m)}^\infty \beta^u \Big(\underline M(t,u)-m\Big)^+ 
	\ + \right. \\
&\qquad\qquad\qquad\qquad \left. + \sum_{u=t}^{\theta_\sigma(m)-1}\beta^u
	\Big(\underline M(t,u)-m\Big)\bigmid\F(\theta_\sigma(m))\right] \\
&= \beta^{\theta_\sigma(m)}D(\theta_\sigma(m))
	+ (1-\beta)\sum_{u=t}^{\theta_\sigma(m)-1}\beta^u\Big(\underline M(t,u)-m\Big),
\end{align*}
making use of \eqref{eq:U-3}, the $\F(u)$-measurability of $\underline M(t,u)$, and the fact that $\underline M(t,u)>m$ for every $u<\sigma(t;m)$.
		
On the other hand, observe that the left-hand side of the identity
$$m\beta^t 
	= m\beta^{\theta_\sigma(m)}+\sum_{u=t}^{\theta_\sigma(m)-1}\beta^u(1-\beta)m$$
does not depend on $\theta$, so the sequence
\begin{equation}\label{eq:U-6}
\left\{\beta^{\theta_\sigma(m)} \Big[V(\theta_\sigma(m);m)-m\Big] 
	+ \sum_{u=t}^{\theta_\sigma(m)-1}
	\beta^u\Big(h(u+1) - (1-\beta)m\Big)\right\}_{\theta=t}^\infty
\end{equation}
is an $\{\F(\theta)\}_{\theta=t}^\infty$-martingale, by \eqref{eq:Z-1} and Proposition \ref{prop:Z-mtg}. It follows that the difference of \eqref{eq:U-5} and \eqref{eq:U-6}, namely,
\begin{multline}\label{eq:U-7}
\left\{\beta^{\theta_\sigma(m)} \Big[V(\theta_\sigma(m); 
	\vphantom{\sum_{u=t}^{\theta_\sigma(m)-1}} \right.
	\underline M(t;\theta_\sigma(m))) - \underline M(t;\theta_\sigma(m))\Big]\ + \\
	\left. + \sum_{u=t}^{\theta_\sigma(m)-1}
	\beta^u \Big(h(u+1)-(1-\beta)\underline M(t,u)\Big)\right\}_{\theta=t}^\infty,
\end{multline}
is an $\{\F(\theta)\}_{\theta=t}^\infty$-martingale, for every $m\ge0$. By \eqref{eq:sigma-0}, $\theta_\sigma(m)\equiv\theta$ when $m=0$; taking $m=0$ in \eqref{eq:U-7} completes the proof.
\end{proof}

	%
	%
	
	\section{The dynamic allocation problem}\label{s:dap}
	
	We now formulate the dynamic allocation problem with general filtrations, {\em with no assumptions about independence or conditional independence whatsoever}, and present a solution in the simple case where the rewards sequences are almost surely decreasing. In particular, using the forward induction characterization of the Gittins index \eqref{eq:gittins-forwards} for each project, we derive a representation of the maximum total expected rewards and show that it is attained by allocation strategies which satisfy the so-called `synchronization identity.' Later, in Section \ref{s:strategies}, we explore the implications of the synchronization identity and how it relates to Gittins indices.
		 
	Let $(\Omega,\F,\P)$ be a complete probability space with a complete multi-parameter filtration $\FF:=\{\F(\tilde s)\}_{\tilde s\in\N_0^d}$, in the sense that $\FF$ satisfies
	\begin{enumerate}[label=(F\arabic*),ref=F\theenumi] \itemsep 1mm
		\item $\F(\tilde s)\subseteq\F(\tilde r)$ for every $\tilde s,\tilde r\in\N_0^d$, 
			$\tilde s\le\tilde r;$ \label{itm:F1}
		\item $\F(\tilde 0)$ contains all the null sets in $\F$; \label{itm:F2}
	\end{enumerate}
	where we put the canonical partial ordering on $\N_0^d$,
	$$\tilde s\le\tilde r\ \iff\ s_i\le r_i,\ \forall\,i=1,\dots,d, \qquad \tilde r,\tilde s\in\N_0^d.$$
	Multi-parameter filtrations in this formulation were first introduced in \cite{Cairoli1975}.
	
	\begin{example}\label{ex:joins}
	A simple example of a multi-parameter filtration is the one constructed as the join of $d$ complete filtrations $\{\F_i(t)\}_{t\in\N_0}$, $i=1,\dots,d$, namely
	$$\F(\tilde s) := \bigvee_{i=1}^d \F_i(s_i) = \sigma\left(\bigcup_{i=1}^d \F_i(s_i)\right),
		\qquad \tilde s = (s_1,\dots,s_d)\in\N_0^d.$$
	\end{example}
	
	We also define two varieties of one-parameter filtrations in which information is only accumulated with respect to a single project, while information regarding the other projects is either complete or frozen. In the former case, we set
	\begin{equation}\label{eq:big-filtr}
		\F^i(t):=\F(\infty,\dots,t,\dots,\infty) 
			= \sigma\left(\bigcup_{\substack{\tilde r\in\N_0^d\\r_i\le t}}\F(\tilde r)\right),
			\qquad t\in\N_0,
	\end{equation}
	where the $t$ appears in the $i$th coordinate; notice that $\F(\tilde s)=\bigcap_{i=1}^d\F^i(s_i)$.
	Similarly, in the latter case, we set
	\begin{equation}\label{eq:small-filtr}
		\F_i(t):=\F(s_1,\dots,s_{i-1},t,s_{i+1}\dots,s_d),\qquad t\ge s_i,
	\end{equation}
	with the convention $\F_i(t):=\F_i(t;\tilde 0)$. In other words, $\F^i(t)$ represents the history of the $i$th project up to time $t$, with the future evolution of all other projects known in advance at time 0, and we set $\FF^i:=\{\F^i(t)\}_{t\in\N_0}$. Analogously, $\F_i(t)$ represents the history of the $i$th project up to time $t$, with this project considered in isolation from all the rest frozen at position $\tilde s$, and we set $\FF_i:=\{\F_i(t)\}_{t\in\N_0}$.
	
	\begin{definition}\label{def:strategy}
	Given any $\tilde s\in\N_0^d$, an {\bf allocation strategy} for $\tilde s$ is an $\N_0^d$-valued sequence $\tilde T:=\{\tilde T(t)\}_{t\in\N_0}$ such that
	\begin{enumerate}
		\item $\tilde T(0) = \tilde s$;
		\item $\tilde T(t+1) = \tilde T(t)+\tilde e_j$ for some $j\in\{1,\dots,d\}$, where $\tilde e_j$ is the unit vector in $\N_0^d$ in the $j$th coordinate; and
		\item $\{\tilde T(t+1) = \tilde T(t)+\tilde e_j,\ \tilde T(t)=\tilde r\}\in\F(\tilde r)$ for every $j=1,\dots,d$.
	\end{enumerate}
	Let $\A(\tilde s)$ denote the collection of allocation strategies for $\tilde s$.

	The first condition requires that $\tilde T$ begin at $\tilde s$, while the second condition implies that $\tilde T$ is a random path in $\N_0^d$ that only takes steps in a positive direction along the integer lattice; the $T_i(t)$'s count the number of times the $i$th project has been engaged up to time $t$. The third condition, also known as the `non-anticipativity' requirement, says that the decision about which project to engage at each time $t$ is made on the basis of the information accumulated from engaging various projects up to time $t$. Notice that condition (2) is equivalent to the requirement that each component $T_i(\cdot)$ be increasing with the constraint $\sum_{i=1}^d(T_i(t)-T_i(0))=t$ for every time $t\in\N_0$.	
	\end{definition}
	
	\begin{definition}
		We call a measurable function $\tilde\nu:\Omega\to\N_0^d$ a {\bf stopping point} of the 
		multi-parameter filtration $\FF$, if $\{\tilde\nu=\tilde r\}\in\F(\tilde r)$ for every 
		$\tilde r\in\N_0^d$. We define also the $\sigma$-algebra generated by a stopping point 
		$\tilde\nu$ to be
		\begin{equation}\label{eq:stopping-point-filtr}
			\F(\tilde\nu) := \{A\in\F\colon A\cap\{\tilde\nu=\tilde r\}\in\F(\tilde r),
			\ \forall\,\tilde r\in\N_0^d\}.
		\end{equation}
		It can be shown that $\tilde\nu$ is $\F(\tilde\nu)$-measurable. This definition generalizes 
		the notion of stopping time from one-parameter filtrations to the multi-parameter setting.
	\end{definition}
	
	\begin{proposition}\label{prop:multiparam-filtr}
	Let $\tilde T\in\A(\tilde s)$ be an allocation strategy. For every $t\in\N_0$, $\tilde T(t)$ is an $\FF$-stopping point, and $T_i(t)$ is an $\FF^i$-stopping time for each $i=1,\dots,d$. Moreover, $\FF(\tilde T):=\{\F(\tilde T(t))\}_{t\in\N_0}$ is a one-parameter filtration, and $\tilde T(\nu)$ is an $\FF$-stopping point for every $\FF(\tilde T)$-stopping time $\nu$.
	\end{proposition}
	\begin{proof}
		By conditions (2) and (3) in Definition \ref{def:strategy}, we obtain
		\begin{equation}\label{eq:multiparam-filtr-1}
			\{\tilde T(t)=\tilde r\} 
				= \bigcup_{j=1}^d\{\tilde T(t+1)=\tilde T(t)+\tilde e_j,\ \tilde T(t)=\tilde r\}
				\in \F(\tilde r)
		\end{equation}
		for every $t\in\N_0$ and $\tilde r\in\N_0^d$ such that $\tilde r\ge\tilde s$. If $\tilde r\not
		\ge\tilde s$ then it follows from condition (2) of Definition \ref{def:strategy} 
		that $\{\tilde T(t)=\tilde r\}=\varnothing\in\F(\tilde r)$. Thus $\tilde T(t)$ is a 
		stopping point. Next, fix $r\in\N_0$. For every $\tilde r\in\N_0^d$ such that 
		$r_i=r$, we have $\F(\tilde r)\subseteq\F^i(r)$, hence
		$$\{T_i(t)=r\} 
			= \bigcup_{\substack{\tilde r\in\N_0^d\\r_i=r}}\{\tilde T(t)=\tilde r\} \in \F^i(r).$$
		Thus $T_i(t)$ is an $\FF^i$-stopping time.
		
		Condition (2) also implies that $\tilde T(t+1)\ge\tilde T(t)$ in the partial ordering
		of $\N_0^d$ for every $t\in\N_0$. Then $\F(\tilde T(t))\subset\F(\tilde T(t+1))$, and thus
		$\{\F(\tilde T(t))\}_{t\in\N_0}$ is a one-parameter filtration.
		
		Finally, let $\nu$ be an $\FF(\tilde T)$-stopping time. For every $\tilde r\in\N_0^d$, we have
		$$\{\tilde T(\nu)=\tilde r\} = \bigcup_{t=0}^\infty\{\nu=t\}\cap\{\tilde T(t)=\tilde r\}.$$
		Since $\{\nu=t\}\in\F(\tilde T(t))$ then \eqref{eq:stopping-point-filtr} implies
		$\{\nu=t\}\cap\{\tilde T(t)=\tilde r\}\in\F(\tilde r)$. We conclude that 
		$\{\tilde T(\nu)=\tilde r\}\in\F(\tilde r)$, which proves that $\tilde T(\nu)$ is an 
		$\FF$-stopping point.
	\end{proof}
	
	\begin{definition}
		A {\bf policy available at $\tilde s\in\N_0^d$} is a pair $(\tilde T,\eta)$, where 
		$\tilde T\in\A(\tilde s)$ is an allocation strategy for $\tilde s$ and 
		$\eta:\Omega\to\N_0\cup\{+\infty\}$ is an $\FF(\tilde T)$-stopping time. We denote by 
		$\p(\tilde s)$ the collection of policies available at $\tilde s$. In other words, a policy 
		is an allocation strategy with an additional option of abandoning all projects at a random 
		time $\eta$, based on the flow of information generated by the strategy.
	\end{definition}
	
	We are now ready to formulate the dynamic allocation problem. Fix the discount factor $\beta\in(0,1)$ and consider $d$ positive, $\FF_i$-predictable rewards sequences $H_i:=\{h_i(t)\}_{t\in\N}$, $i=1,\dots,d$, one for each project. We assume that each $h_i(t)$ is uniformly bounded from above by $K(1-\beta)$ for some fixed $K>0$, so that the maximum discounted reward that can accumulate over time for any given project is $\sum_{t=0}^\infty\beta^tK(1-\beta) = K$. 
	
	\begin{definition}
		For retirement reward $M\in[0,\infty)$ and $\tilde s\in\N_0^d$, the {\bf dynamic allocation 
		problem} is to compute the value
		\begin{equation}\label{eq:dap}
			\Phi(\tilde s;M) := \esssup_{(\tilde T,\eta)\in\p(\tilde s)}
			\,\E[\reward(\tilde T,\eta;M)\mid\F(\tilde s)],
		\end{equation}
		where
		\begin{equation}\label{eq:dap-reward}
			\reward(\tilde T,\eta;M) := \sum_{i=1}^d\sum_{t=0}^{\eta-1} \beta^t h_i(T_i(t+1))
				\Big(T_i(t+1)-T_i(t)\Big) + M\beta^\eta
		\end{equation}
		is the total discounted reward accumulated by following the allocation strategy $\tilde T$ with
		a retirement reward $M$ collected at the $\FF(\tilde T)$-stopping time $\eta$. As a 
		consequence of bounding $h_i(t)$ within the interval $[0,K(1-\beta)]$, the conditional 
		expectation $\E[\reward(\tilde T,\eta;M)\mid\F(\tilde s)]$ is well-defined and bounded.
	\end{definition}
	
	In solving the dynamic allocation problem, we consider the family of optimal stopping problems
	$$V_i(t;m) := \esssup_{\tau\in\s_i(t)}\,\E\left[\sum_{u=t}^{\tau-1}\beta^{u-t} h_i(u+1)
		+ m\beta^{\tau-t} \bigmid \F_i(t) \right],\qquad i=1,\dots,d,$$
	indexed by $m\ge0$, associated with the previously-established filtrations $\FF_i$ and reward
	sequences $H_i$. Here, $\s_i(t)$ is the collection of $\FF_i$-stopping times taking values in 
	$\{t,t+1,\dots\}\cup\{+\infty\}$. Each such optimal stopping problem has an associated Gittins 
	index sequence $\{M_i(t)\}_{t\in\N_0}$ and its lower envelope $\underline M_i(t,\theta)$ for 
	$t,\theta\in\N_0$, $t\le\theta$. Notice that bounding each $h_i(t)$ within $[0,K(1-\beta)]$ 
	implies, via \eqref{eq:osp}, that $V_i(t;m)\le K$ for each $i=1,\dots,d$, hence
	\begin{equation}\label{eq:gittins-upper}
		M_i(t) \le K,\qquad \text{a.s.}
	\end{equation}
	as a consequence of \eqref{eq:gittins}, for every $i=1,\dots,d$ and $t\in\N_0$. All the results from 
	Sections \ref{s:osp} and \ref{s:gittins} hold for these new objects.
	
	To capture the combined effects of all $d$ projects, we further introduce, for every $m\ge0$ and $\tilde s\in\N_0^d$, the `total operational time'
	\begin{equation}\label{eq:tau}
		\tau(m;\tilde s) := \sum_{i=1}^d(\sigma_i(s_i;m)-s_i),\qquad\tau(m) := \tau(m;\tilde 0),
	\end{equation}
	and its right-inverse
	\begin{equation}\label{eq:N}
		N(t;\tilde s) := \inf\{m\ge0\colon \tau(m;\tilde s)\le t, \text{ a.s.}\},\qquad 
		N(t) := N(t;\tilde 0); \qquad t\in\N_0.
	\end{equation}
	Here, $\tau(m;\tilde s)$ represents the earliest calendar time at which it will be optimal to abandon all projects and collect the retirement reward $m$, when starting with position $\tilde s=(s_1,\dots,s_d)$ at calendar time $t=\sum_{i=1}^d s_i$. Analogously, $N(t;\tilde s)$ represents the minimum reward that must be offered to abandon all projects immediately, at the current position $\tilde s$ and current calendar time $t=\sum_{i=1}^d s_i$. In other words, $N(t;\tilde s)$ is the current `equitable surrender value for the entire collection of projects.'
	
	The analogues of \eqref{eq:sigma-equivs} and \eqref{eq:sigma-M}, namely
	\begin{equation}\label{eq:tau-equivs}
		\tau(m;\tilde s) > t\ \iff\ m < N(t;\tilde s),\qquad\text{a.s.}
	\end{equation}
	and
	\begin{equation}\label{eq:tau-N}
		\tau(m;\tilde s) = \inf\{\theta\ge t\colon N(\theta;\tilde s)\le m,\ \text{a.s.}\},
	\end{equation}
	now hold for every $m\ge0$, $\tilde s=(s_1,\dots,s_d)\in\N_0^d$, and $t=\sum_{i=1}^d s_i$. We have also the analogue of \eqref{eq:sigma-jumps},
	\begin{equation}\label{eq:tau-jumps}
		\tau(N(t;\tilde s);\tilde s)\le t<\tau(N(t;\tilde s)-;\tilde s),\qquad\text{a.s.}
	\end{equation}
	as well as the following analogue of Proposition \ref{prop:M-constant}.
	
	\begin{proposition}\label{prop:N-constant}
		If $\tau(\,\cdot\,;\tilde s)$ has a jump at $m\ge0$ for fixed $\tilde s\in\N_0^d$, then 
		the mapping $\theta\mapsto N(\theta;\tilde s)$ is almost surely constant on 
		$\{\tau(m;\tilde s),\dots,\tau(m-;\tilde s)-1\}$, and
		\begin{equation}\label{eq:N-constant}
			m = N(\tau(m;\tilde s);\tilde s) = N(\tau(m;\tilde s)+1;\tilde s)
				= \dots = N(\tau(m-;\tilde s)-1;\tilde s),\quad\text{a.s.}
		\end{equation}
	\end{proposition}
	
	It will be useful to utilize the following representation of $N(t;\tilde s)$, as a more explicit version of Proposition \ref{prop:N-constant}. For ease of notation, we assume $\tilde s=\tilde 0$.
	
	Notice that $\tau(\cdot)$ is a decreasing, integer-valued step function, as a consequence of Lemma \ref{lem:sigma-props}, with $\tau(0) = \infty$ and $\tau(m)\to0$ as $m\to\infty$. Then there exist pathwise finitely many times $0 = t_0 < t_1 < \dots < t_\ell = \infty$ with corresponding values $m_1 > \dots > m_\ell \ge 0$ such that $\tau(\cdot)$ has jumps precisely at $m_1,\dots,m_\ell$, with $\tau(m_j) = t_{j-1}$ for $j=1,\dots,\ell$. It follows from \eqref{eq:tau-equivs} that $N(\cdot)$ is a bounded, decreasing, discrete-time process satisfying
	\begin{equation}\label{eq:N-step}
		N(t) = \sum_{j=1}^\ell m_j\1\{t_{j-1}\le t<t_j\},\qquad t\in\N_0,
	\end{equation}
	as well as
	\begin{equation}\label{eq:tau-step}
		t\in\{t_{j-1},\dots,t_j-1\}\ \implies\ \tau(N(t))=t_{j-1}\,\text{ and }\,\tau(N(t)-)=t_j.
	\end{equation}
	We remark that the times $t_1,\dots,t_\ell$ and the values $m_1,\dots,m_\ell$ are of course
	random---determined pathwise for the function $\tau(\cdot)$---as is $\ell$ itself.
	
	\subsection*{Dynamic allocation with decreasing rewards}
	
	We consider the simple yet indispensable case where the rewards sequences $\widehat H_i := \{\widehat h_i(t)\}_{t\in\N}$ are almost surely decreasing, while remaining positive and $\FF_i$-predictable. We wish to compute the value of the dynamic allocation problem
	\begin{equation}\label{eq:dap-decr}
		\widehat\Phi(\tilde s;M) := \esssup_{(\tilde T,\eta)\in\p(\tilde s)}
			\,\E[\widehat\reward(\tilde T,\eta;M)\mid\F(\tilde s)], \qquad\tilde s\in\N_0^d,
	\end{equation}
	where
	\begin{equation}\label{eq:dap-decr-reward}
		\widehat\reward(\tilde T,\eta;M) := \sum_{i=1}^d \sum_{t=0}^{\eta-1} \beta^t
			\widehat h_i(T_i(t+1))\Big(T_i(t+1)-T_i(t)\Big) + M\beta^\eta.
	\end{equation}
	We emphasize that what follows requires absolutely no assumptions on the dependency structure of the multi-parameter filtration $\FF$.
	
	As above, we have for each project a corresponding optimal stopping problem, for which the objects and results found in Sections \ref{s:osp} and \ref{s:gittins} stand. It turns out that the lower envelopes of the Gittins indices of each project relate directly to the values of their associated rewards sequences, via the forward induction characterization of the Gittins index.
		
	\begin{lemma}\label{lem:decr-reward-rep}
		For each $i=1,\dots,d$ and $t\ge s_i$, we have the representation
		\begin{equation}\label{eq:decr-reward-rep}
			\widehat h_i(t+1) = (1-\beta)\underline M_i(s_i,t),\qquad\text{a.s.}
		\end{equation}
	\end{lemma}
	\begin{proof}
		For every $t\ge s_i$ and stopping time $\nu\in\s(t+1)$, we have by the decrease and 
		predictability of $\widehat h_i(\cdot)$ that
		$$\E\left[\sum_{u=t}^{\nu-1}\beta^u\widehat h_i(u+1)\bigmid\F_i(t)\right]
			\le \widehat h_i(t+1)\,\E\left[\sum_{u=t}^{\nu-1}\beta^u\bigmid\F_i(t)\right],\qquad\text{a.s.}$$
		It follows from \eqref{eq:gittins-forwards} that $(1-\beta)M_i(t)\le \widehat h_i(t+1)$ almost surely. 
		But taking $\nu\equiv t+1$ in the essential supremum of \eqref{eq:gittins-forwards} also yields 
		$(1-\beta)M_i(t)\ge \widehat h_i(t+1)$ almost surely. Now $M_i(\cdot)$ is decreasing, and since 
		$\underline M_i(s_i,\cdot)$ is the lower envelope of $M_i(\cdot)$, then in fact 
		$M_i(t)=\underline M_i(s_i,t)$, for every $t\ge s_i$.
	\end{proof}
	
	We are now ready to provide a solution to the dynamic allocation problem \eqref{eq:dap-decr}, conditioned on the existence of an allocation strategy $\tilde T\in\A(\tilde s)$ that satisfies the `synchronization identity,' i.e.,
	\begin{equation}\label{eq:sync}
		\sum_{i=1}^d \Big[(T_i(t)-s_i)\wedge(\sigma_i(s_i;m)-s_i)\Big] = t\wedge\tau(m;\tilde s),
		\qquad \forall\,t\in\N_0,\,m>0.
	\end{equation}
	We will construct an explicit strategy satisfying \eqref{eq:sync} in Section \ref{s:strategies}.
	
	\begin{theorem}\label{thm:dap-decr-soln}
		Suppose that the multi-parameter filtration $\FF$ satisfies conditions \eqref{itm:F1} and \eqref{itm:F2}. Then the value of the dynamic allocation problem \eqref{eq:dap-decr} is given by
	\begin{multline}\label{eq:dap-decr-V}
		\widehat\Phi(\tilde s;M) = M + \E\left[\int_M^\infty \left(1-\beta^{\tau(m;\tilde s)}\right) \dd m
			\bigmid \F(\tilde s)\right] \\
		= (1-\beta)\,\E\left[\sum_{t=0}^\infty \beta^t (M\vee N(t;\tilde s))
			\bigmid \F(\tilde s)\right], \quad\text{a.s.}
	\end{multline}
	Moreover, the supremum in \eqref{eq:dap-decr} is attained by $(\tilde T,\eta)\in\p(\tilde s)$ if and only if $\tilde T\in\A(\tilde s)$ satisfies the synchronization identity \eqref{eq:sync} for every $m\ge M$, and either $\eta=\tau(M;\tilde s)$ or $\eta=\tau(M;\tilde s)+1$.
	\end{theorem}
	\begin{proof}
		Without loss of generality, we assume that $\tilde s=\tilde 0$. For every $t\in\N_0$, $m\ge M$, and 
		$\tilde T\in\A(\tilde 0)$, we define
		$$A_i(t;m,\tilde T):=T_i(t)\wedge\sigma_i(m),
			\qquad A(t;m,\tilde T):=\sum_{i=1}^dA_i(t;m,\tilde T)\le t\wedge\tau(m).$$
		Notice that $A(t;m,\tilde T)=t\wedge\tau(m)$ holds if and only if $\tilde T$ satisfies the 
		synchronization identity \eqref{eq:sync} for $m$. Also notice that $A_i(0;m,\tilde T)=0$ and that
		\begin{equation}\label{eq:dap-decr-1}
			\1\{\sigma_i(m) > T_i(t)\} \Big(T_i(t+1)-T_i(t)\Big)
				= A_i(t+1;m,\tilde T)-A_i(t;m,\tilde T).
		\end{equation}
		We compute the following, using \eqref{eq:sigma-equivs}, \eqref{eq:decr-reward-rep}, and \eqref{eq:dap-decr-1}:
\begin{align*}
	\widehat\reward(\tilde T&,\eta;M) \\
	&= \sum_{i=1}^d \sum_{t=0}^{\eta-1} \beta^t(1-\beta) \Big(T_i(t+1)-T_i(t)\Big) \,\cdot \\
	&\qquad\qquad \cdot \left(M\wedge\underline M_i(T_i(t))+\int_M^\infty 
		\1\{m < \underline M_i(T_i(t))\}\,\dd m\right) + M\beta^\eta \\
	&\le \sum_{i=1}^d \int_M^\infty \sum_{t=0}^{\eta-1} \beta^t(1-\beta) 
		\Big(A_i(t+1;m,\tilde T) - A_i(t;m,\tilde T)\Big)\,\dd m + M \\
	&= \int_M^\infty \sum_{t=0}^{\eta-1} \beta^t(1-\beta) 
		\Big(A(t+1;m,\tilde T) - A(t;m,\tilde T)\Big)\,\dd m + M \\
	&\le \int_M^\infty\left[(1-\beta)\beta^{\eta-1}(\eta\wedge\tau(m)) 
		+ (1-\beta)^2\sum_{t=0}^{\eta-2} \beta^t \Big((t+1)\wedge\tau(m)\Big)\right]\,\dd m + M.
\end{align*}
The first inequality becomes equality if and only if $M\le\underline M_i(T_i(t))$ for each $i=1,\dots,d$ and $0\le t<\eta$. This is equivalent to $t\le\tau(M)$ for $t<\eta$ by \eqref{eq:sigma-equivs}, hence we require $\eta\le\tau(M)+1$. The second inequality follows from expanding the telescoping sum and becomes equality if and only if $\tilde T$ satisfies the synchronization identity. By splitting the sum into $t<\tau(m)$ and $t\ge\tau(m)$ observing that the integrand is increasing in $\eta$, one observes that the integrand is bounded from above by $1-\beta^{\tau(m)}$. But it is easy to check that this upper bound is attained when $\eta\ge\tau(m)$ for all $m\ge M$. Since $m\mapsto\tau(m)$ is decreasing by Lemma \ref{lem:sigma-props}, it suffices that $\eta\ge\tau(M)$. This proves the first equality in \eqref{eq:dap-decr-V}. The second equality of \eqref{eq:dap-decr-V} follows from \eqref{eq:tau-equivs}, since
\begin{multline*}
	(1-\beta) \sum_{t=0}^\infty \beta^t (M\vee N(t))
	= (1-\beta) \sum_{t=0}^\infty \beta^t\left(M + \int_M^\infty\1\{t<\tau(m)\}\,\dd m\right) \\
	= M + \int_M^\infty \left(1-\beta^{\tau(m)}\right) \dd m.
\end{multline*}

Finally, to show that $\tau(M)$ is an $\FF(\tilde T)$-stopping time, observe that \eqref{eq:sync} implies
$$\{\tau(M)\le t\} = \bigcap_{i=1}^d \{\sigma_i(M)\le T_i(t)\} 
	= \{\tilde\sigma(m) := (\sigma_1(m),\dots,\sigma_d(m))\le\tilde T(t)\}.$$
Recall that $\tilde T(t)$ is an $\FF$-stopping point, from Proposition \ref{prop:multiparam-filtr}, and notice that
$$\{\tilde\sigma(m)\le\tilde r\} = \bigcap_{i=1}^d\{\underline M_i(r_i)\le m\}
	\in \bigcap_{i=1}^d \F^i(r_i) = \F(\tilde r), \qquad\tilde r\in\N_0^d$$
by \eqref{eq:sigma-equivs} and the definition \eqref{eq:big-filtr}. Thus
$$\{\tau(M)\le t\}\cap\{\tilde T(t)=\tilde r\} = \{\tilde\sigma(m)\le\tilde r\}\in\F(\tilde r)$$
and so $\{\tau(M)\le t\}\in\F(\tilde T(t))$ by definition, for each $t\in\N_0$.
\end{proof}
	
In words, Theorem \ref{thm:dap-decr-soln} states that the satisfying the synchronization identity is a necessary and sufficient condition for an allocation strategy to be optimal, and that one should retire as soon as every rewards sequence, or, in light of \eqref{eq:decr-reward-rep}, the Gittins index of every project, falls below the retirement reward $M$.

	%
	%
	
	\section{Synchronization and index-type strategies}\label{s:strategies}
	Before continuing further, we elucidate additional characterizations of the synchronization identity \eqref{eq:sync}, which will simplify later proofs and reveal new representations of the value of the dynamic allocation problem. We also introduce allocation strategies of `index type'---namely, strategies which compel the investigator to engage, at each moment in time, a project with the maximum current Gittins index---and show via the `excursion theory' of \cite{Kaspi1998} that index-type strategies are a subset of those satisfying the synchronization identity. Finally, we construct an explicit strategy of index type that therefore satisfies the synchronization identity and completes the proof of Theorem \ref{thm:dap-decr-soln}.

	First, let us introduce the random fields
	\begin{equation}\label{eq:M-cal}
		\M(\tilde s):=\max_{1\le i\le d} M_i(s_i),
	\end{equation}
	and
	\begin{equation}\label{eq:M-cal-lower}
		\underline\M(\tilde s,\tilde r):=\max_{1\le i\le d}\underline M_i(s_i,r_i),
		\qquad \underline\M(\tilde r):=\underline\M(\tilde 0,\tilde r)
	\end{equation}
	for every $\tilde s,\tilde r\in\N_0^d$ such that $\tilde s\le\tilde r$, recalling the notation
	of \eqref{eq:gittins} and \eqref{eq:gittins-envelope}.
	
	\begin{theorem}[Synchronization and dual optimality]\label{thm:sync-props}
		Let $\tilde T\in\A(\tilde s)$ be an allocation strategy. The following are equivalent:
		\begin{enumerate}[label=(\alph*)]
			\item $\tilde T$ satisfies the synchronization identity \eqref{eq:sync};
			\item for each $m\ge0$,
				\begin{equation}\label{eq:sync-1}
					T_i(\tau(m;\tilde s)) = \sigma_i(s_i;m),\qquad i=1,\dots,d;
				\end{equation}
			\item for each $t\in\N_0$,
				\begin{equation}\label{eq:sync-2}
					\sigma_i(s_i;N(t;\tilde s)) \le T_i(t) \le \sigma_i(s_i;N(t;\tilde s)-),
						\qquad i=1,\dots,d;
				\end{equation}
			\item $\tilde T$ satisfies the `dual optimality property,' i.e.,
				\begin{equation}\label{eq:sync-dual}
					\max_{1\le i\le d}\underline M_i(s_i,T_i(t)) = \underline\M(\tilde s,\tilde T(t)) 
						= N(t;\tilde s),\qquad\forall\,t\in\N_0;
				\end{equation}
			\item $\tilde T$ is of the `lower envelope index type,' i.e.,
				\begin{equation}\label{eq:lower-index-type}
					\underline M_i(s_i,T_i(t)) = \underline\M(\tilde s,\tilde T(t))
						\quad\text{on}\quad \{\tilde T(t+1) = \tilde T(t) + \tilde e_i\}
				\end{equation}
				for each $i=1,\dots,d$ and $t\in\N_0$.
		\end{enumerate}
	\end{theorem}
	
	Of particular importance here is the synchronization identity \eqref{eq:sync}, which will feature heavily in Section \ref{s:f4}. In words, a strategy that satisfies the synchronization identity will not engage any given project beyond the optimal stopping time associated with that project, no matter how large or small the exit reward, until {\em every} project has been engaged up to its optimal stopping time. Compare this with \eqref{eq:sync-1}, which posits that the strategy $\tilde T$ partitions the total operational time $\tau(m;\tilde s)$ optimally across all projects, for any choice of retirement reward and starting position. 
	
	On the other hand, we shall see below in the proof of Theorem \ref{thm:sync-props} that 
	\begin{equation}\label{eq:dual-bound}
		\underline\M(\tilde s,\tilde T(t)) \ge N(t;\tilde s)
		\tag{\ref*{eq:sync-dual}$'$}
	\end{equation}
	holds for {\em every} allocation strategy $\tilde T\in\A(\tilde s)$; the dual optimality property \eqref{eq:sync-dual} is then saying that strategies satisfying the synchronization identity are not only optimal for the maximization problem in dynamic allocation, as in Section \ref{s:dap}, but they also solve a related {\em minimization} problem, expressed by \eqref{eq:sync-dual} and \eqref{eq:dual-bound}.
	
	Notice also that \eqref{eq:sync} is equivalent to the similar statement
	\begin{equation}\label{eq:sync-rcts}
		\sum_{i=1}^d\Big[T_i(t)-s_i)\wedge(\sigma_i(s_i;m-)-s_i)\Big] = t\wedge\tau(m-;\tilde s),
			\qquad\forall\,t\in\N_0,\,m>0.
	\end{equation}
	by the right-continuity of $\sigma_i(s_i;\,\cdot\,)$ and $\tau(\,\cdot\,;\tilde s)$.
	
	\begin{proof}[Proof of Theorem \ref{thm:sync-props}]
		We will proceed by proving (a)$\implies$(b)$\implies$(c)$\implies$(a), then show 
		(a)$\implies$(d)$\implies$(b) and finally (a)$\iff$(e). Throughout the proof, we assume
		for ease of notation that $\tilde s=\tilde 0$, as the general case follows by translation.
		
		\vspace{2mm}
		(a)$\implies$(b). 
		Taking $t=\tau(m)$ in \eqref{eq:sync} yields $\sum_{i=1}^d(T_i(\tau(m))
		\wedge\sigma_i(m))=\tau(m)$, hence $T_i(\tau(m))=\sigma_i(m)$ for each $i=1,\dots,d$ and every
		$m\ge0$, by condition (2) of Definition \ref{def:strategy}.
		
		\vspace{2mm}
		(b)$\implies$(c). 
		It follows from \eqref{eq:tau-jumps}, \eqref{eq:sync-1}, and the increase of $T_i(\cdot)$ that,
		for each $i=1,\dots,d$,
		$$\sigma_i(N(t))=(T_i\circ\tau)(N(t))\le T_i(t)\le(T_i\circ\tau)(N(t)-)=\sigma_i(N(t)-).$$
		
		\vspace{2mm}
		(c)$\implies$(a). Fix $t\in\N_0$ and $m\ge0$. If $m>N(t)$ then \eqref{eq:sync-2} and the
		decrease of $\sigma_i(\cdot)$ yield 
		$$T_i(t)\ge\sigma_i(N(t))\ge\sigma_i(m)$$
		Now \eqref{eq:tau-jumps} and the decrease of $\tau(\cdot)$ imply 
		\begin{equation}\label{eq:sync-3}
			t\ge\tau(N(t))\ge\tau(m),
		\tag{\ref*{eq:sync-rcts}$'$}
		\end{equation}
		yielding \eqref{eq:sync}.
		
		On the other hand, if $m\le N(t)$, we have analogously that
		\begin{equation}\label{eq:sync-4}
			T_i(t)\le\sigma(N(t)-)\le\sigma_i(m-)\quad\text{and}\quad t<\tau(N(t)-)\le\tau(m-),
		\tag{\ref*{eq:sync-rcts}$''$}
		\end{equation}
		yielding \eqref{eq:sync-rcts}. But \eqref{eq:sync-rcts} is equivalent to 
		\eqref{eq:sync}.
	
		\vspace{2mm}
		(a)$\implies$(d).
		Let $\tilde T$ be an arbitrary allocation strategy and let $t\in\N_0$ and $m\ge0$ be given. 
		It holds that $\underline\M(\tilde T(t))\le m$ if and only if $\underline M_i(T_i(t))\le m$ 
		for each $i=1,\dots,d$, which in turn holds if and only if $\sigma_i(m)\le T_i(t)$ for each 
		$i=1,\dots,d$, by \eqref{eq:sigma-equivs}. Summing over each $i$ yields $\tau(m)\le t$, which is 
		equivalent to the inequality $N(t)\le m$ by \eqref{eq:tau-equivs}. Thus 
		\begin{equation}\label{eq:dual-1}
			\underline\M(\tilde T(t))\le m\ \implies\ N(t)\le m,\qquad\forall\,t\in\N_0,\,m\ge0,
		\end{equation}
		from which \eqref{eq:dual-bound} follows. Now the reverse inequality 
		holds if and only if the converse of \eqref{eq:dual-1} also holds, which requires the implication
		\begin{equation}\label{eq:dual-2}
			\tau(m)\le t\ \implies\ \sigma_i(m)\le T_i(t),
				\qquad\forall\,i=1,\dots,d,\ t\in\N_0,\,m\ge0.
		\end{equation}
		But this holds immediately for every $\tilde T$ that satisfies \eqref{eq:sync}.
		
		\vspace{2mm}
		(d)$\implies$(b).
		If \eqref{eq:sync-dual} holds then \eqref{eq:dual-2} must also hold.
		By taking $t=\tau(m)$, we obtain $\sigma_i(m)\le T_i(\tau(m))$ for each
		$i=1,\dots,d$ and every $m\ge0$. But summing over $i$ yields
		$$\tau(m)=\sum_{i=1}^d\sigma_i(m)\le\sum_{i=1}^dT_i(\tau(m))=\tau(m),$$
		hence $\sigma_i(m)=T_i(\tau(m))$ for each $i=1,\dots,d$.
		
		\vspace{2mm}
		(a)$\implies$(e).
		It suffices to show that, for each $i=1,\dots,d$ and $t\in\N_0$, we have
		$$\underline M_i(T_i(t)) < \underline\M(\tilde T(t))\ \implies\ T_i(t) = T_i(t+1).$$
		We have by property (d) that $\underline\M(\tilde T(t))=N(t)$, so we assume 
		$\underline M_i(T_i(t)) < N(t)$ for given $i=1,\dots,d$ and $t\in\N_0$.
		
		Suppose $\sigma_i(\cdot)$ has a jump at $N(t)$. Then $T_i(t)\in\{\sigma_i(N(t)),\dots,
		\sigma_i(N(t)-)\}$, by \eqref{eq:sync-2}. If $T_i(t)\neq\sigma_i(N(t)-)$ then
		\eqref{eq:M-constant} implies that $\underline M_i(T_i(t)) = N(t)$, contradicting our 
		initial assumption. Thus, $T_i(t)=\sigma_i(N(t)-)$ holds. In particular, we have by 
		\eqref{eq:sync-1} that $T_i(t) = T_i(\tau(N(t)-))$, where $t<\tau(N(t)-)$ by 
		\eqref{eq:tau-jumps}. The increase of $T_i(\cdot)$ thus implies that $T_i(t)=T_i(t+1)$, as desired.
		
		On the other hand, suppose $\sigma_i(\cdot)$ is continuous at $N(t)$. Then 
		\eqref{eq:sync-2} implies that $T_i(t)=\sigma_i(N(t))=\sigma_i(N(t)-)$. Notice that, 
		by \eqref{eq:sync-1}, we have
		$$T_i(\tau(N(t)-)) = \sigma_i(N(t)-) = T_i(t),$$
		whereas $\tau(N(t)-)>t$ by \eqref{eq:tau-jumps}. The increase of $T_i(\cdot)$ again implies that 
		$T_i(t)=T_i(t+1)$.
		
		\vspace{2mm}
		(e)$\implies$(a).
		Fix $m>0$ and proceed by induction on $t\in\N_0$. For $t=0$, we have $T_i(0)=0$ for each 
		$i=1,\dots,d$, and synchronization is trivially satisfied. Now suppose synchronization is 
		satisfied for some $t\in\N_0$. We consider two cases.
				
		{\it Case 1}: $\sum_{i=1}^d(T_i(t)\wedge\sigma_i(m))=\tau(m)\le t$. This holds if and only if 
		$\sigma_i(m)\le T_i(t)$ holds for each $i=1,\dots,d$. The increase of $T_i(\cdot)$ then implies
		$\sum_{i=1}^d(T_i(t+1)\wedge\sigma_i(m))=\tau(m)\le t+1$, as desired.
		
		{\it Case 2}: $\sum_{i=1}^d(T_i(t)\wedge\sigma_i(m))=t\le\tau(m)$. This holds if and only if
		$\sigma_i(m)\ge T_i(t)$ holds for each $i=1,\dots,d$. If equality holds for each $i$ then
		we fall into Case 1. Therefore, suppose there exists a non-empty set $J\subseteq\{1,\dots,d\}$
		such that $\sigma_j(m)>T_j(t)$ holds if $j\in J$ and $\sigma_j(m)=T_j(t)$ holds if $j\notin J$.
		By \eqref{eq:sigma-equivs}, $\underline M_j(T_j(t))>m$ holds if $j\in J$ whereas 
		$\underline M_j(T_j(t))\le m$ holds if $j\notin J$. It follows from the lower envelope 
		index-type property \eqref{eq:lower-index-type} that $\tilde T(t+1) = \tilde T(t)+e_{j^*}$ for 
		some $j^*\in J$. Then $\sigma_i(m)\ge T_i(t)=T_i(t+1)$ for each $i\neq j^*$ while 
		$\sigma_{j^*}(m)>T_{j^*}(t)$ implies $\sigma_{j^*}(m)\ge T_{j^*}(t)+1=T_{j^*}(t+1)$. As a 
		result, $\sigma_i(m)\ge T_i(t+1)$ holds for each $i=1,\dots,d$, and synchronization is 
		satisfied for $t+1$.
	\end{proof}
	
	\begin{remark}
		The synchronization identity \eqref{eq:sync} can also be understood as allowing the 
		`minimum' and `summation' operations to commute, when it comes to the pairs of random 
		times $T_i(t)$ and $\sigma_i(s_i;m)$. One interesting connection here is with the 
		max-plus algebra (also known as the tropical semiring) of \cite{Baccelli1992}; see also 
		\cite{ElKaroui2008}. Alternatively, the explicit relationship between synchronization and index-type strategies, as we show below in Proposition \ref{prop:excursion}, may help elucidate the role of the Gittins index in certain scheduling problems of queueing theory; see, for example, \cite{Scully2018,Scully2025,Ramakrishna2026}.
	\end{remark}
	
	\begin{remark}\label{rmk:myopic}
	In light of Theorem \ref{thm:sync-props}, particularly \eqref{eq:sync-dual}, we immediately observe from \eqref{eq:dap-decr-V} that the value of the dynamic allocation problem \eqref{eq:dap-decr} with decreasing rewards sequences has the additional representation
	\begin{equation}\label{eq:myopic}
		\widehat\Phi(\tilde s;M) = (1-\beta)\,\E\left[\sum_{t=0}^\infty \beta^t (\underline\M(\tilde s,\tilde T^*(t))\vee M) \bigmid \F(\tilde s)\right],
	\end{equation}
	for any $\tilde T^*\in\A(\tilde s)$ satisfying the synchronization identity \eqref{eq:sync}. The equivalence of the synchronization identity and the lower envelope index property \eqref{eq:lower-index-type}, in conjunction with \eqref{eq:decr-reward-rep}, shows that optimal strategies for the decreasing rewards case are exactly those which engage, at each moment in time, the project that offers the maximum immediate reward. Such strategies are sometimes referred to as `myopic' in the sense that one does not need to consider the future at all when deciding which project to engage. A similar proof that myopic strategies are optimal in this special setting also appeared in \cite{Bank2007}.
	\end{remark}
	
	\subsection*{Index-type strategies and excursions}
	
	We now introduce allocation strategies of index type and connect them to the synchronization identity by means of the `excursion-theoretic' approach to dynamic allocation first developed in \cite{Kaspi1994,Kaspi1998}.
	
	\begin{definition}
	An allocation strategy $\tilde T\in\A(\tilde s)$ is of {\bf index type} if
	\begin{equation}\label{eq:index-type}
		M_i(T_i(t)) = \M(\tilde T(t)) \qquad\text{on}\quad \{\tilde T(t+1)=\tilde T(t)+\tilde e_i\},
	\end{equation}
	recalling the notation of \eqref{eq:M-cal}. In other words, an index-type strategy requires the investigator to engage a project only if it has maximal current Gittins index.
	\end{definition}
	
	Compare \eqref{eq:index-type} against the weaker \eqref{eq:lower-index-type}, which requires the investigator to only engage projects whose current \emph{minimum Gittins index to date} is maximal amongst all projects. This observation is the basis for the following connection between index-type strategies and strategies which satisfy the synchronization identity. The proof can be found in the \hyperref[appendix]{Appendix}.
	
	\begin{proposition}\label{prop:excursion}
	An allocation strategy $\tilde T\in\A(\tilde s)$ is of index type \eqref{eq:index-type} if and only if it satisfies the synchronization identity \eqref{eq:sync} and the `excursion property'
	\begin{equation}\label{eq:excursion}
		\{M_i(T_i(t)) > \underline M_i(T_i(t))\} \subseteq \{\tilde T(t+1) = \tilde T(t) + \tilde e_i\}.
	\end{equation}
	\end{proposition}
	
	We term \eqref{eq:excursion} the `excursion property' in the tradition of \cite{Kaspi1994}. Morally speaking, when a project is engaged, its Gittins index may either go above its lower envelope, or it may fall to a new `all-time low.' Periods of time in which the Gittins index exceeds its lower envelope are referred to as `excursions' away from the lower envelope; excursions come to an end when the Gittins index drops to a new all-time low and again coincides with its lower envelope. As such, the excursion property \eqref{eq:excursion} is satisfied by an allocation strategy if the investigator exclusively engages the project that is in the middle of an excursion period. When all projects' Gittins indices coincide with their lower envelopes, such as at time $t=0$, the investigator has the option, \emph{prima facie}, to engage any one of them; the synchronization identity and its equivalent \eqref{eq:lower-index-type} then obliges the investigator to choose a project with maximal current Gittins index and to pursue that project until the end of its excursion.
	
	For later use in the proof of Theorem \ref{thm:main}, let us expand briefly on the idea of excursions. Let $\FF$ be a complete multi-parameter filtration. Recall for each $i=1,\dots,d$ that $\s_i(t)$ denotes the collection of $\FF_i$-stopping times taking values in $\{t,t+1,\dots\}\cup\{+\infty\}$. For some fixed $t\in\N_0$, we introduce the stopping times $\varepsilon_i(k;t)\in\s_i(t)$, $k\in\N_0$, defined by $\varepsilon_i(0;t)\equiv t$ and 
	\begin{equation}\label{eq:eps}
	\varepsilon_i(k+1;t) 
		:= \inf\{\theta>\varepsilon_i(k;t)\colon M_i(\theta)\le M_i(\varepsilon_i(k;t))\},
		\qquad k\in\N_0.
	\end{equation}
	Clearly 
	$$t=\varepsilon_i(0;t) < \varepsilon_i(1;t) < \cdots, \quad\text{and}\quad
		\lim_{k\to\infty}\varepsilon_i(k;t) = +\infty, \qquad\text{a.s.}$$
	so, in particular, we have
	\begin{equation}\label{eq:eps-set}
	\{\varepsilon_i(k;t)\colon k\in\N_0\} = \{\theta\in\N_0\colon M_i(\theta)=\underline M_i(t,\theta)\}, 
		\qquad\text{a.s.},
	\end{equation}
	In other words, the stopping time $\varepsilon_i(k;t)$ denotes the beginning of the $k$th excursion of the Gittins index, from time $t$ onwards, for the $i$th project.
	
	Additionally, for each $\tilde s,\tilde k\in\N_0^d$, let 
	\begin{equation}\label{eq:eps-tilde}
	\tilde\varepsilon(\tilde k;\tilde s) 
		:= (\varepsilon_1(k_1;s_1),\dots,\varepsilon_d(k_d;s_d))\in\N_0^d.
	\end{equation}
	For any allocation strategy $\tilde T\in\A(\tilde s)$, define $\rho(k;\tilde T)$, $k\in\N_0$ by $\rho(0;\tilde T)\equiv0$ and
	\begin{equation}\label{eq:rho}
	\rho(k+1;\tilde T) := \inf\{\theta > \rho(k;\tilde T) 
		\colon \tilde T(\theta) = \tilde\varepsilon(\tilde k;\tilde s) \text{ for some } \tilde k\in\N_0^d\},
	\end{equation}
	which denote the calendar times at which, in the course of pursuing the allocation strategy $\tilde T$, every project is about to begin its next excursion. The next result, whose proof is in the \hyperref[appendix]{Appendix}, states that the random times $\rho(k;\tilde T)$ are accessible to an investigator pursuing strategy $\tilde T$.
	
	\begin{proposition}\label{prop:rho-stopping}
	For every $\tilde k,\tilde s\in\N_0^d$, $\tilde\varepsilon(\tilde k;\tilde s)$ is an $\FF$-stopping point, and for every $k\in\N_0$, $\rho(k;\tilde T)$ is an $\FF(\tilde T)$-stopping time.
	\end{proposition}
	
	In light of Proposition \ref{prop:rho-stopping}, we define the random, possibly infinite `excursion intervals' of an allocation strategy $\tilde T\in\A(\tilde s)$ by
	\begin{equation}\label{eq:excursion-interval}
	E(k;\tilde T) := \{\rho(k;\tilde T),\dots,\rho(k+1;\tilde T)-1\}, \qquad k\in\N_0.
	\end{equation}
	Our terminology is motivated by the following crucial lemma.

	\begin{lemma}\label{lem:rho-excursion}
	If an allocation strategy $\tilde T\in\A(\tilde s)$ satisfies the excursion property \eqref{eq:excursion}, then for every $k\in\N_0$, there exists pathwise a unique $j\equiv j(k)\in\{1,\dots,d\}$ such that $\tilde T(t+1)=\tilde T(t)+\tilde e_j$ for every $t\in E(k;\tilde T)$.
	\end{lemma}
	\begin{proof}
	Let us work pathwise, i.e., for fixed $\omega\in\Omega$. Fix $k\in\N_0$. At time $t=\rho(k;\tilde T)$, we have $M_i(T_i(\rho(k;\tilde T)))=\underline M_i(T_i(\rho(k;\tilde T)))$ for each $i=1,\dots,d$, by \eqref{eq:eps-set}, \eqref{eq:rho}. Set $j(k)$ to be the project engaged at time $\rho(k;\tilde T)$, i.e., the unique $j$ for which $\tilde T(\rho(k;\tilde T)+1) = \tilde T(\rho(k;\tilde T)) + \tilde e_j$.
	
	For our inductive hypothesis, let $t\in E(k;\tilde T)$ with $t<\rho(k+1;\tilde T)-1$, and suppose that $\tilde T(t+1)=\tilde T(t)+\tilde e_j$ and that $M_i(t)=\underline M_i(t)$ for each $i\neq j$. Since $\rho(k;\tilde T) < t < \rho(k+1;\tilde T)$, then there exists at least one project for which the Gittins index does not coincide with its lower envelope. This must be the $j$th project by assumption, i.e., $M_j(t) > \underline M_j(t)$. The excursion property \eqref{eq:excursion} implies that the inductive hypothesis also holds for $t+1$. It follows that the strategy engages the $j$th project exclusively during the excursion interval $E(k;\tilde T)$.
	\end{proof}
	
	

\subsection*{An explicit index-type strategy}
	
In the remainder of this section, we construct a strategy $\tilde T^*\in\A(\tilde 0)$ satisfying the synchronization identity \eqref{eq:sync}, via an approach strongly suggested by the inequalities \eqref{eq:sync-2}. Similar constructions appear in \cite{Mandelbaum1987, ElKaroui1994, ElKaroui1997}.
	
Let $y_0(t):=\tau(N(t))$, recalling the notation of \eqref{eq:tau} and \eqref{eq:N}, and for each $i=1,\dots,d$, define
\begin{equation}\label{eq:y-traverse}
y_i(t) := \tau(N(t)) + \sum_{j=1}^{i}\Delta\sigma_j(N(t)) 
	= \sum_{j=1}^{i}\sigma_j(N(t)-) + \sum_{j=i+1}^d\sigma_j(N(t)), \qquad t\in\N,
\end{equation}
where $\Delta\sigma_i(m):=\sigma_i(m-)-\sigma_i(m)$. Recall that $\tau(N(t))\le t<\tau(N(t)-)$ from \eqref{eq:sync-3} and \eqref{eq:sync-4}, and we have
\begin{equation}\label{eq:y-partition}
\tau(N(t)) = y_0(t) \le y_1(t) \le \dots \le y_d(t) = \tau(N(t)-).
\end{equation}
Thus, there exists a unique $k\equiv k(t)\in\{1,\dots,d\}$ for which $t\in[y_{k-1}(t),y_k(t))$. Now for each $i=1,\dots,d$, we simply set $T^*_i(0)=0$ and
\begin{equation}\label{eq:y-strategy}
T^*_i(t) = \left. \begin{cases}
	\sigma_i(N(t)-) &\text{if }\, i=1,\dots,k-1 \\
	\sigma_i(N(t)) + t - y_{i-1}(t) &\text{if }\, i=k \\
	\sigma_i(N(t)) &\text{if }\, i=k+1,\dots,d
\end{cases}\right\},
\qquad t\in\N.
\end{equation}

\begin{proposition}\label{prop:y-strategy}
The process $\tilde T^*(t)=(T^*_1(t),\dots,T^*_d(t))$, as defined in \eqref{eq:y-strategy}, is an allocation strategy of index type \eqref{eq:index-type}.
\end{proposition}
\begin{proof}
Conditions (1) and (3) of Definition \ref{def:strategy} are trivially satisfied. It is clear, furthermore, that $\sum_{i=1}^dT^*_i(t)=t$ for every $t\in\N_0$. To show that each $T^*_i(\cdot)$ is increasing, fix $i=1,\dots,d$ and recall the notation of \eqref{eq:N-step}. Suppose for some $j'=1,\dots,\ell$ that $t\in\{t_{j'-1},\dots,t_{j'}-1\}$ . We wish to show that $T^*_i(t)\le T^*_i(t+1)$.
		
Suppose $t+1<t_{j'}$. Then $N(t)=N(t+1)$ and so $y_j(t)=y_j(t+1)$ for each $j=0,\dots,d$. Additionally, it holds that $k(t)\le k(t+1)$, so it suffices to show that
\begin{equation}\label{eq:y-strategy-1}
\sigma_i(N(t)) \le \sigma_i(N(t))+t-y_{k(t)-1}(t) \le \sigma_i(N(t)-).
\end{equation}
The first inequality follows from the definition of $k(t)$. Similarly, we have
$$t\le y_{k(t)}(t)=\sigma_i(N(t)-)-\sigma_i(N(t))+y_{k(t)-1}(t),$$
which is equivalent to the second inequality.

Suppose on the other hand that $t+1=t_{j'}$. But \eqref{eq:y-strategy-1} and the decrease of both $\sigma_i(\cdot)$ and $N(\cdot)$ imply
$$T^*_i(t) \le \sigma_i(N(t)-) \le \sigma_i(N(t+1)) \le T^*_i(t+1),$$
hence $\tilde T$ is an allocation strategy.

Finally, suppose $T_i^*(t)=\sigma_i(N(t))$. Since $\sigma_i(\cdot)$ is integer-valued, it is a right-continuous step function. Thus, $\sigma_i(\cdot)$ has a jump at the point $r = \inf\{m\ge0\colon \sigma_i(m)=\sigma_i(N(t))\}$, which implies that $\underline M_i(\sigma_i(r))=r$, by \eqref{eq:M-sigma}. But \eqref{eq:M-sigma} implies $M(\sigma_i(r))=r$ as well, yielding $M_i(T_i^*(t)) = \underline M_i(T_i^*(t))$. If instead $T_i^*(t)=\sigma_i(N(t)-)$ then we may repeat the same argument by instead setting $r = \inf\{m\ge0 \colon \sigma_i(m) = \sigma_i(N(t)-)\}$. It follows that, in the notation of \eqref{eq:y-strategy}, if $i\neq k$, then the excursion property \eqref{eq:excursion} holds vacuously. Since \eqref{eq:y-strategy} guarantees that the $k$th project is engaged, in fact \eqref{eq:excursion} holds in general. Additionally, $\tilde T^*$ satisfies \eqref{eq:sync} as a consequence of \eqref{eq:y-strategy-1} and Theorem \ref{thm:sync-props}. Thus $\tilde T^*$ is of index type, by Proposition \ref{prop:excursion}.
\end{proof}
	
	Roughly speaking, the allocation strategy \eqref{eq:y-strategy} is a time change via $\sigma_i(N(\cdot))$. The `flat' parts of $N(\cdot)$ correspond to the jumps of $\tau(\cdot)$, which are in turn comprised of the jumps of each $\sigma_i(\cdot)$. Immediately after every jump of $\tau(\cdot)$, our strategy has the form $T^*_i(t)=\sigma_i(N(t))$. To ensure that $T^*_i(t)$ lies between $\sigma_i(N(t))$ and $\sigma_i(N(t)-)$, the investigator should engage each project, {\em one at a time} in `lexicographic order,' until the distance between $\sigma_i(N(t))$ and $\sigma_i(N(t)-)$ has been `traversed.' That this traversal occurs one project at a time ensures that $\tilde T^*$ satisfies the excursion property. We remark that in the language of excursion theory, this strategy may be described as an `index-priority strategy' which adheres to a `static priority scheme' $(1,\dots,d)$ \cite{Kaspi1998}.
	
\begin{remark}
While the construction \eqref{eq:y-strategy} seems quite complex at first glance, in practice it is quite straightforward to implement since it is of index type.  At time $t=0$, one computes the Gittins indices of each project and selects a project with maximal index to engage. One continues to engage this same project until its Gittins index is no longer maximal amongst all projects, at which point one switches once more to a project whose Gittins index is maximal. In case of ties, \eqref{eq:y-strategy} is designed to select the project lowest in lexicographic order, although of course one can use any manner of tie-breaking. In general, the benefit of index policies from a computational perspective is that at each moment in time, one only needs to compute $d$ indices to make optimal decisions, rather than searching the entire space of paths through $\N_0^d$, thus resolving the so-called `curse of dimensionality.' On the other hand, computing the Gittins index can be a non-trivial task, although algorithms have been analyzed in the finite-state Markov setting \cite{Varaiya1985}. For an excellent overview of practical implementations of Gittins index policies, see \cite{Scully2025}.
\end{remark}

\section{General multi-parameter martingales}\label{s:multiparam}
	
We set now to solve the dynamic allocation problem \eqref{eq:dap} in the general setting, without additional assumptions on the rewards sequences as in \eqref{eq:dap-decr}. In contrast to previous works, the filtrations $\FF_1,\dots,\FF_d$ associated with each project are not necessarily independent. As before, we work on a complete probability space $(\Omega,\F,P)$ with a multi-parameter filtration $\FF:=\{\F(\tilde s)\}_{\tilde s\in\N_0^d}$. This filtration is assumed to have the properties \eqref{itm:F1} and \eqref{itm:F2}, as introduced in Section \ref{s:dap}, as well as the following additional property of \cite{Cairoli1975}:
\begin{enumerate}[start=4,label=(F\arabic*),ref=F\theenumi]\itemsep 1mm
\item \label{itm:F4} $\F(\tilde s),\F(\tilde r)$ are conditionally independent given 
	$\F(\tilde s\wedge\tilde r)$, meaning that 
	\begin{equation}\label{eq:cond-indep}
	\P(A\cap B\mid\F(\tilde s\wedge\tilde r))
		= \P(A\mid\F(\tilde s\wedge\tilde r))\,\P(B\mid\F(\tilde s\wedge\tilde r))
	\end{equation}
holds for every $A\in\F(\tilde s)$ and $B\in\F(\tilde r)$.
\end{enumerate}
Notice that the (F3) right-continuity condition in \cite{Cairoli1975} is omitted here because we are in a discrete time setting. The celebrated \eqref{itm:F4} condition is standard in multi-parameter theory, and constitutes a relaxation of the projects' independence by allowing them to mutually depend on a `shared history.' We use the unnatural ennumeration, in passing over (F3), to remain consistent with the well-established appellation `(F4)' used for the condition spelled out in \eqref{eq:cond-indep}.
	
\begin{proposition}[Walsh \cite{Walsh1986}, p.\ 349]\label{eq:F4-equivs}
Let $\FF$ be a multi-parameter filtration. The following are equivalent:
\begin{enumerate}[label=(\alph*)] \itemsep 1mm
\item $\FF$ satisfies condition \eqref{itm:F4};
\item for every $\tilde s,\tilde r\in\N_0^d$ and bounded $\F(\tilde r)$-measurable function $\xi$,
\begin{equation}\label{eq:F4-1}
\E[\xi\mid\F(\tilde s)]=\E[\xi\mid\F(\tilde s\wedge\tilde r)];
\end{equation}
\item for every $\tilde s,\tilde r\in\N_0^d$ and bounded $\F$-measurable function $\xi$,
\begin{equation}\label{eq:F4-2}
\E[\,\E[\xi\mid\F(\tilde s)]\mid\F(\tilde r)]
	= \E[\,\E[\xi\mid\F(\tilde r)]\mid\F(\tilde s)].
\end{equation}
\end{enumerate}
If any of the above hold then 
\begin{equation}\label{eq:F4-cap}
\F(\tilde s)\cap\F(\tilde r) = \F(\tilde s\wedge\tilde r).
\end{equation}
holds for every $\tilde s,\tilde r\in\N_0^d$.
\end{proposition}
Proofs for these claims can be found in \cite{Cairoli1975} and \cite{Kallenberg1997}, while proofs for the claims made in the following Examples \ref{ex:F4-indep} and \ref{ex:F4-sheet} can be found in the \hyperref[appendix]{Appendix}.
	 
\begin{example}\label{ex:F4-indep}
The multi-parameter filtration constructed as the join of $d$ one-parameter filtrations as in Example \ref{ex:joins}, where the one-parameter filtrations are mutually independent, satisfies condition \eqref{itm:F4}.
\end{example}
	
\begin{example}\label{ex:F4-sheet}
We provide an example, due to \cite{Cairoli1975}, of a multi-parameter filtration that satisfies \eqref{itm:F4} but cannot be expressed as a join of independent one-parameter filtrations. Let $\B$ denote the power set $\sigma$-algebra on $\N_0^d$ and let $\{\eta(\tilde s)\}_{\tilde s\in\N_0^d}$ be a family of i.i.d.\ standard Gaussian random variables. Then the so-called `white noise process' on $\N_0^d$ is the family of random variables $\{W(S)\}_{S\in\B,\,|S|<\infty}$, indexed by finite sets in $\B$, defined by $W(S):=\sum_{\tilde s\in S}\eta(\tilde s)$. This family satisfies
\begin{enumerate}\itemsep 1mm
\item $\E[W(S)]=0$, $\E[W(S)^2]=|S|$; and
\item if $S\cap R = \varnothing$ then $W(S)$ and $W(R)$ are independent with $W(S)+W(R)=W(S\cup R)$ almost surely;
\end{enumerate}
and has covariance matrix given by
\begin{equation}\label{eq:sheet-covar}
\E[W(S)W(R)] = \sum_{\substack{\tilde s\in S,\,\tilde r\in R\\\tilde s\neq\tilde r}}
	\E[\eta(\tilde s)]\E[\eta(\tilde r)] 
	+ \sum_{\tilde s\in S\cap R}\E[\eta(\tilde s)]^2 = |S\cap R|.
\end{equation}
The process $\{B(\tilde s)\}_{\tilde s\in\N_0^d}$, defined as $B(\tilde s) := W(\{1,\dots,s_1\}\times\dots\times\{1,\dots,s_d\})$, is the `Brownian sheet' on $\N_0^d$, and the filtration it generates,
\begin{equation}\label{eq:sheet-filtr}
\F(\tilde r) := \sigma\left(B(\tilde s)\colon \tilde s\le\tilde r,\ \tilde s\in\N_0^d\right),
	\qquad\tilde r\in\N_0^d,
\end{equation}
satisfies \eqref{itm:F4}, but cannot be written as the join of independent one-parameter filtrations as a consequence of \eqref{eq:sheet-covar}.
\end{example}

\begin{example}
In light of the previous example, consider the dynamic allocation problem defined on the filtration generated by the 2-dimensional Brownian sheet $\{B(\tilde s)\}_{\tilde s\in\N_0^2}$, in which the rewards processes are given by $h_1(t+1)=B(t,0)^2\wedge K$ and $h_2(t+1)=B(0,t)^2 \wedge K$ for $t\in\N_0$ and some fixed $K>0$. These are respectively $\FF_1$- and $\FF_2$-predictable, bounded, and {\em not independent}, since they both depend on $\eta(0,0)$. Of course, one can construct, in this vein, more involved dependencies between the two rewards sequence.
\end{example}

	Recall the one-parameter filtrations $\FF_i$ and $\FF^i$ defined in \eqref{eq:big-filtr} and \eqref{eq:small-filtr}. The condition \eqref{itm:F4} provides a link between the `small' filtration $\FF_i$ and the `large' filtration $\FF^i$, via the following `enlargement of filtrations' property.
	
	\begin{proposition}\label{prop:enlarge}
		Every $\FF_i$-supermartingale is also an $\FF^i$-supermartingale.
	\end{proposition}
	\begin{proof}
		It suffices to show, for every $s,t\in\N_0$, $s\le t$, and every bounded, $\F_i(t)$-measurable
		random variable $\xi$, that $\E[\xi\mid\F_i(s)] = \E[\xi\mid\F^i(s)]$. Notice by 
		\eqref{eq:F4-cap} that $\F_i(t)\cap\F^i(s)=\F_i(s)$, so \eqref{eq:F4-1} yields
		\begin{equation*}
			\E[\xi\mid\F_i(s)] = \E[\xi\mid\F_i(t)\cap\F^i(s)] = \E[\xi\mid\F^i(s)]
		\end{equation*}
		as desired.
	\end{proof}
	
	\begin{remark}\label{rmk:enlarge}
		As a consequence of the previous proposition, we observe that the results of Sections 
		\ref{s:osp} and \ref{s:gittins} hold with respect to either 
		filtration, $\FF_i$ or $\FF^i$.
	\end{remark}
	
%
	
	We now introduce several new objects in our multi-parameter setting. For the remainder of this section, we assume $\FF$ satisfies conditions \eqref{itm:F1}, \eqref{itm:F2}, and \eqref{itm:F4}.
	
	\begin{definition}
		Let $\mathbf X:\N_0^d\times\Omega\to\R$ be a random field, with $X(\cdot)$ denoting a 
		sample field. We say that $\mathbf X$ is an $\FF${\bf-submartingale} if
		\begin{enumerate}\itemsep 1mm
			\item $\E|X(\tilde s)|<\infty$ for every $\tilde s\in\N_0^d$;
			\item $\mathbf X$ is $\FF$-adapted, i.e., $X(\tilde s)$ is $\F(\tilde s)$-measurable for
				every $\tilde s\in\N_0^d$; and
			\item $\E[X(\tilde r)\mid\F(\tilde s)]\ge X(\tilde s)$ almost surely for
				every $\tilde s\le \tilde r$.
		\end{enumerate}
		We say $\mathbf X$ is an $\FF${\bf-supermartingale} if $-\mathbf X$ is an 
		$\FF$-submartingale, and we say $\mathbf X$ is an $\FF${\bf-martingale} if it is both 
		an $\FF$-submartingale and an $\FF$-supermartingale.
	\end{definition}
	
	\begin{proposition}[Walsh \cite{Walsh1986}, Proposition 2.1, p.\ 349]\label{prop:mp-mtg}
		For a random field $\mathbf X:\N_0^d\times\Omega\to\R$, the following are equivalent:
		\begin{enumerate}[label=(\alph*)]\itemsep 1mm
			\item $\mathbf X$ is an $\FF$-supermartingale;
			\item $\mathbf X$ is $\FF$-adapted, integrable, and 
				$\E[X(\tilde r)\mid\F(\tilde s)]\le X(\tilde s\wedge\tilde r)$ a.s., for every
				$\tilde s,\tilde r\in\N_0^d$;
			\item the process $\{X(\tilde s+t\tilde e_i)\}_{t\in\N_0}$ is a 
				$\{\F^i(s_i+t)\}_{t\in\N_0}$-supermartingale for every $\tilde s\in\N_0^d$ and for 
				each $i=1,\dots,d$.
		\end{enumerate}
	\end{proposition}
	
	The following result, whose proof can be found in the \hyperref[appendix]{Appendix}, says that following an $\FF$-supermartingale along the path of an allocation strategy is itself a one-parameter supermartingale. We remark that this fact is subtler than it may seem at first glance, as the optional sampling theorem used in the usual proof in a single dimension fails to hold for multi-parameter filtrations of dimension $d\ge2$ \cite{Mandelbaum1981}. In fact, it is precisely sequences of stopping points generated by allocation strategies, in the sense of Proposition \ref{prop:multiparam-filtr}, that enjoy an optional sampling-type property. For details on the issues that arise in the multi-dimensional setting, see \cite{Mandelbaum1981,Walsh1981} as well as \cite[Thm.~3]{Mandelbaum1986}.
	
\begin{proposition}[Mandelbaum--Vanderbei \cite{Mandelbaum1981}, p.~255]\label{prop:mp-T-supermtg}
Let $\tilde T\in\A(\tilde s)$ be an allocation strategy and let $\mathbf X$ be an $\FF$-supermartingale. Then $\mathbf X(\tilde T) := \{X(\tilde T(t))\}_{t\in\N_0}$ is an $\FF(\tilde T)$-supermartingale.
\end{proposition}
\begin{proof}
Fix arbitrary $t\in\N_0$, $\tilde r\in\N_0^d$, and $\Lambda\in\F(\tilde T(t))$. Since $\FF(\tilde T)$ is a filtration by Proposition \ref{prop:multiparam-filtr}, it suffices to show $\E[X(\tilde T(t+1))\1_{\Lambda\cap\{\tilde T(t)=\tilde r\}}]\le\E[X(\tilde T(t))\1_{\Lambda\cap\{\tilde T(t)=\tilde r\}}]$ by definition of conditional expectation. We have
\begin{multline*}
\E\left[X(\tilde T(t+1))\1_{\Lambda\cap\{\tilde T(t)=\tilde r\}}\right] 
= \sum_{i=1}^d \E\left[X(\tilde r+\tilde e_i)\1_{\Lambda\cap\{\tilde T(t)=\tilde r,\,\tilde T(t+1)=\tilde T(t)+\tilde e_i\}}\right] \\
\le \sum_{i=1}^d \E\left[X(\tilde r)\1_{\Lambda\cap\{\tilde T(t)=\tilde r,\,\tilde T(t+1)=\tilde T(t)+\tilde e_i\}}\right]
\le \E\left[X(\tilde T(t))\1_{\Lambda\cap\{\tilde T(t)=\tilde r\}}\right],
\end{multline*}
where the inequality follows from the supermartingality of $\bX$ with the fact that $\Lambda\cap\{\tilde T(t)=\tilde r,\,\tilde T(t+1)=\tilde T(t)+\tilde e_i\}\in\F(\tilde r)$ by definition of allocation strategies.
\end{proof}

	To conclude this section, we prove a result that will be needed in Section \ref{s:f4}.
		
	\begin{lemma}\label{lem:K-mtg}
		Let $\tilde T\in\A(\tilde s)$ be an allocation strategy. For each $i=1,\dots,d$, suppose
		$Q_i=\{Q_i(u)\}_{u=s_i}^\infty$ is an $\{\F_i(u)\}_{u=s_i}^\infty$-(super)martingale. Then the process
		$$\sum_{i=1}^d \sum_{u=0}^{t-1} \beta^{u-T_i(u)} \Big(Q_i(T_i(u+1))-Q_i(T_i(u))\Big),
			\qquad t\in\N_0$$
		is an $\FF(\tilde T)$-(super)martingale.
	\end{lemma}
	\begin{proof}
		We prove the lemma for supermartingales; the same arguments hold for the martingale case. We may assume $\tilde s=\tilde 0$. For each $i=1,\dots,d$, the process
		$$U_i(t) := \sum_{u=0}^{t-1}\beta^{-u}\Big(Q_i(u+1)-Q_i(u)\Big),\qquad t\in\N_0$$
		is an $\FF_i$-supermartingale because it is an integrable martingale transform of 
		$Q_i$ by the deterministic process $\{\beta^{-u}\}_{u\in\N_0}$. In particular, it is an
		$\FF^i$-supermartingale by the enlargement of filtrations property, Proposition
		\ref{prop:enlarge}. We obtain then from Proposition \ref{prop:mp-mtg} (the 
		(c)$\implies$(a) implication) that the random field
		$$U(\tilde r):=\sum_{i=1}^dU_i(r_i),\qquad \tilde r\in\N_0^d,$$
		is an $\FF$-supermartingale. Now Proposition \ref{prop:mp-T-supermtg} implies that
		\begin{align*}
			U(\tilde T(t)) 
			&= \sum_{i=1}^d\sum_{u=0}^{t-1}\beta^{-T_i(u)}\Big(Q_i(T_i(u+1))-Q_i(T_i(u))\Big),
				\qquad t\in\N_0
		\end{align*}
		is an $\FF(\tilde T)$-supermartingale. Taking the martingale transform of $U(\tilde T)$ by  $\{\beta^u\}_{u\in\N_0}$ yields the desired $\FF(\tilde T)$-supermartingale.
	\end{proof}

	\section{Optimality in the \eqref{itm:F4} setting}\label{s:f4}
	
	We now solve the general dynamic allocation problem \eqref{eq:dap}, with respect to a multi-parameter filtration $\FF$ satisfying only the conditions \eqref{itm:F1}, \eqref{itm:F2}, and \eqref{itm:F4}. As before, each project, viewed in isolation, corresponds to an optimal stopping problem, for which the objects and results of Sections \ref{s:osp} and \ref{s:gittins} stand.
	
	The main result of this work, Theorem \ref{thm:main}, says that the general dynamic allocation problem \eqref{eq:dap} has the same value as the dynamic allocation problem \eqref{eq:dap-decr}, in which the rewards sequences are {\em almost surely decreasing}, and that all index-type strategies, which are optimal in the decreasing rewards setting, are {\em also optimal in the general case}. This result will require almost all of the machinery developed in Sections \ref{s:strategies} and \ref{s:multiparam} and, in particular, relies crucially on the excursion property \eqref{eq:excursion} described in Proposition \ref{prop:excursion}.
	
	Versions of Theorem \ref{thm:main} in which the filtrations $\FF_1,\dots,\FF_d$ are assumed to be independent have appeared in \cite{Mandelbaum1986,ElKaroui1993,Kaspi1998}; in all three of these previous works, independence of filtrations is crucial to the proofs. One should also notice that Proposition \ref{prop:ost-decr-rep} is a special case of Theorem \ref{thm:main}.
	
%
	
	\begin{theorem}\label{thm:main}
		Assume that the multi-parameter filtration $\FF$ satisfies the conditions \eqref{itm:F1}, \eqref{itm:F2}, and \eqref{itm:F4}. Let $\underline M_i(\cdot,\cdot)$, $i=1,\dots,d$ denote the lower envelopes of the corresponding Gittins index sequences $M_i(\cdot)$ defined with respect to the projects of the general dynamic allocation problem \eqref{eq:dap}. 
		
		Then, for every $\tilde s\in\N_0^d$, the value of the dynamic allocation problem \eqref{eq:dap} is almost surely equal to the value of the dynamic allocation problem \eqref{eq:dap-decr} with decreasing rewards given by $\widehat h_i(t+1) := (1-\beta)\underline M_i(s_i,t)$, i.e.,
		\begin{multline}\label{eq:dap-V}
			\Phi(\tilde s;M) = \widehat\Phi(\tilde s;M)
				= M + \E\left[\int_M^\infty \left(1-\beta^{\tau(m;\tilde s)}\right) \dd m 
					\bigmid \F(\tilde s)\right] \\
				= (1-\beta)\,\E\left[\sum_{t=0}^\infty \beta^t(M\vee N(t;\tilde s))
					\bigmid \F(\tilde s)\right],
				\qquad\text{a.s.}
		\end{multline}
		
		Moreover, any policy which follows an allocation strategy $\tilde T^*\in\A(\tilde s)$ of index type \eqref{eq:index-type}, and stops at $\tau(M;\tilde s)$, is optimal for the dynamic allocation problem \eqref{eq:dap}, i.e.,
		\begin{equation}
			\Phi(\tilde s;M)
				= \E[\reward(\tilde T^*) \mid \F(\tilde s)]
				\ge \E[\reward(\tilde T) \mid \F(\tilde s)],
				\quad\text{a.s.},\qquad\forall\,\tilde T\in\A(\tilde s),
		\end{equation}
		and yields the representation
		\begin{equation}\label{eq:main-rep}
		\Phi(\tilde s;M) = (1-\beta)\,\E\left[\sum_{t=0}^\infty \beta^t 
			\Big(M\vee\underline\M(\tilde s,\tilde T^*(t))\Big) \bigmid \F(\tilde s)\right], 
			\qquad\text{a.s.}
		\end{equation}
	\end{theorem}
	
	\begin{remark}
		We must be careful when considering the objects $\tau(\,\cdot\,;\tilde s)$ and 
		$N(\,\cdot\,;\tilde s)$ in \eqref{eq:dap-V}, since they depend on the underlying filtration and rewards sequences---recall their definitions \eqref{eq:tau} and \eqref{eq:N}. In \eqref{eq:dap-V}, they refer to the optimal time to abandon all projects and to the equitable surrender value for all projects, respectively, for the dynamic allocation problem \eqref{eq:dap} with filtration $\FF$ satisfying conditions \eqref{itm:F1}, \eqref{itm:F2}, and \eqref{itm:F4} and general rewards sequences.
	\end{remark}
	
	We prove the special case where no retirement reward is given, i.e., $M=0$, from which the general case can be obtained immediately by incorporating the computations from the proof of Theorem \ref{thm:dap-decr-soln}. We rely on three key facts, namely, that if $\tilde T^*\in\A(\tilde s)$ is an index-type allocation strategy, the following hold:
	\begin{enumerate}[label=(T\arabic*),ref=T\arabic*]\itemsep 1mm	
	\item \label{itm:t1} 
	in the notation of Proposition \ref{prop:U-mtg}, the process
		\begin{equation}\label{eq:t1}
			K(t) := \sum_{i=1}^d \sum_{u=0}^{t-1} \beta^{u-T^*_i(u)}
				\Big(U_i(T_i^*(u+1))-U_i(T^*_i(u))\Big), \qquad t\in\N_0
		\end{equation}
		is a uniformly integrable $\FF(\tilde T^*)$-martingale;
		
	\item \label{itm:t2}
	$K(\infty) := \lim_{t\to\infty}K(t)$ and $\Lambda(\infty) := \lim_{t\to\infty}\Lambda(t)$ 
		exist and are almost surely equal, where
		\begin{multline}
			\Lambda(t) := \sum_{i=1}^d \sum_{u=0}^{t-1} \beta^u 
				\Big[h_i(T_i^*(u+1)) - (1-\beta)\underline M_i(s_i,T_i^*(u))\Big]
				\Big(T_i^*(u+1)-T_i^*(u)\Big), \\ \qquad t\in\N_0,
		\end{multline}
		and the series $\Lambda(\infty)$ is absolutely convergent almost surely; and
	
	\item \label{itm:t3}
	for any allocation strategy $\tilde T\in\A(\tilde s)$, the process
	$$Q(t;\tilde T) := \beta^t\widehat\Phi(\tilde T(t))
		+ \sum_{i=1}^d\sum_{u=0}^{t-1}\beta^u h_i(T_i(u+1))\Big(T_i(u+1)-T_i(u)\Big),
		\qquad t\in\N_0$$
	is a non-negative, bounded $\FF(\tilde T)$-supermartingale.
	\end{enumerate}
	Indeed, \eqref{itm:t1} and \eqref{itm:t2} imply almost surely that 
	$$\E[\Lambda(\infty)\mid\F(\tilde s)] = \E[K(\infty)\mid\F(\tilde s)] = K(0) = 0.$$
	Due to the absolute convergence of $\Lambda(\infty)$, we may write, in the notation of 
	\eqref{eq:dap-reward} and \eqref{eq:dap-decr-reward}, $\Lambda(\infty) = \reward(\tilde T^*) -
	\widehat\reward(\tilde T^*)$. Thus, we have
	\begin{equation}\label{eq:main-1}
		\E[\reward(\tilde T^*)\mid\F(\tilde s)] = \E[\widehat\reward(\tilde T^*)\mid\F(\tilde s)] 
			= \widehat\Phi(\tilde s), \qquad\text{a.s.,}
	\end{equation}
	since $\tilde T^*$ is optimal for the problem \eqref{eq:dap-decr} by Proposition
	\ref{prop:excursion} and Theorem \ref{thm:dap-decr-soln}.
	
	Now for any allocation strategy $\tilde T\in\A(\tilde s)$, let $\{Q(t;\tilde T)\}_{t\in\N_0}$ 
	be the supermartingale of \eqref{itm:t3}. Note that $Q(\infty;\tilde T) := 
	\lim_{t\to\infty}Q(t;\tilde T)$ exists and equals $\reward(\tilde T)$, since $\widehat\Phi(\cdot)$ is 
	bounded. Therefore, we obtain
	$$\widehat\Phi(\tilde s) = \widehat\Phi(\tilde T(0)) = Q(0;\tilde T) \ge \E[Q(\infty;\tilde T)\mid\F(\tilde s)]
		= \E[\reward(\tilde T)\mid\F(\tilde s)], \qquad\text{a.s.},$$
	which in conjunction with \eqref{eq:main-1} yields
	$$\E[\reward(\tilde T^*)\mid\F(\tilde s)] \ge \E[\reward(\tilde T)\mid\F(\tilde s)],
		\quad\text{a.s.},\qquad\forall\,\tilde T\in\A(\tilde s).$$
	By definition, $\tilde T^*$ attains the supremum in \eqref{eq:dap}, and $\Phi(\tilde s)=
	\E[\reward(\tilde T^*)\mid\F(\tilde s)]=\widehat\Phi(\tilde s)$. The representation \eqref{eq:main-rep} follows from \eqref{eq:sync-dual} of Theorem \ref{thm:sync-props}, and the theorem is proved.
	
	It now remains to prove \eqref{itm:t1}, \eqref{itm:t2}, and \eqref{itm:t3}.
	
	\begin{proof}[Proof of \eqref{itm:t2}]
	For each $i=1,\dots,d$, let
	\begin{equation}\label{eq:y-diff}
		Y_i(t) := V_i(t;\underline M_i(T_i^*(t))) - \underline M_i(T_i^*(t)), \qquad t\in\N_0.
	\end{equation}
	In the notation of Proposition \ref{prop:U-mtg}, we may write
	\begin{multline*}
	K(t) = \sum_{i=1}^d \sum_{u=0}^{t-1} \beta^{u-T_i^*(u)} 
		\Big[\beta^{T_i^*(u+1)}Y_i(u+1) - \beta^{T_i^*(u)}Y_i(u)\ + \\
	+ \beta^{T_i^*(u)}\Big(h_i(T_i^*(u+1)) 
		- (1-\beta)\underline M_i(s_i,T_i^*(u))\Big)\Big]
		\Big(T_i^*(u+1)-T_i^*(u)\Big),
	\end{multline*}
	yielding
	\begin{equation}\label{eq:main-2}
	W(t) := K(t)-\Lambda(t) = \sum_{i=1}^d\sum_{u=0}^{t-1}\beta^{u-T_i^*(u)}
		\Big(\beta^{T_i^*(u+1)}Y_i(u+1)-\beta^{T_i^*(u)}Y_i(u)\Big).
	\end{equation}
	Observe that $Y_i(t)=0$ whenever $M_i(T_i^*(t)) = \underline M_i(T_i^*(t))$ by virtue of \eqref{eq:gittins}. Then by \eqref{eq:eps-set}, \eqref{eq:rho}, and Lemma \ref{lem:rho-excursion}, the summand in \eqref{eq:main-2} is non-zero for exactly a single project at each time $u$, and consequently, the almost sure limit $W(\infty) := \lim_{t\to\infty}W(t)$ exists and equals the sum
	\begin{equation}\label{eq:main-3}
	W(\infty) = \sum_{i=1}^d \sum_{k=0}^\infty \sum_{u=\rho(k;\tilde T^*)}^{\rho(k+1;\tilde T^*)-1}	
		\beta^{u-T_i^*(u)}\Big(\beta^{T_i^*(u+1)}Y_i(u+1) - \beta^{T_i^*(u)}Y_i(u)\Big),
	\end{equation}
	in the notation of \eqref{eq:N-step} and \eqref{eq:rho}. Since $|\beta^{-T_i^*(u)}(\beta^{T_i^*(u+1)}Y_i(u+1) - \beta^{T_i^*(u)}Y_i(u))|$ is bounded by $2K$, the sum converges absolutely.
	
	In light of Lemma \ref{lem:rho-excursion}, let $j\equiv j(k)$ be the exclusive project engaged between times $\rho(k;\tilde T^*)$ and $\rho(k+1;\tilde T^*)$. The summand in \eqref{eq:main-3} is zero for every $i\neq j$, and
	\begin{multline}\label{eq:main-4}
	\sum_{u=\rho(k;\tilde T^*)}^{\rho(k+1;\tilde T^*)-1}
		\beta^{u-T_j^*(u)}\Big(\beta^{T_j^*(u+1)}Y_j(u+1)-\beta^{T_j^*(u)}Y_j(u)\Big) \\
	= \beta^{\rho(k+1;\tilde T^*)}Y_j(\rho(k+1;\tilde T^*)) 
		- \beta^{\rho(k;\tilde T^*)}Y_j(\rho(k;\tilde T^*)).
	\end{multline}
	Both terms on the right-hand side are zero since $\rho(k;\tilde T^*),\rho(k+1;\tilde T^*)$ are by definition times where the Gittins indices of every project coincides with their lower envelopes, so we obtain from  \eqref{eq:main-3} that $W(\infty)=0$. Assuming that \eqref{itm:t1} holds, we have the existence of $K(\infty)$, since uniformly integrable martingales converge almost surely. It follows that $\Lambda(\infty)$ exists and 
		$K(\infty)=\Lambda(\infty)$. Moreover, it is clear from the definition that $\Lambda(\cdot)$ 
		is uniformly bounded from above by the random variable
		$$\sum_{i=1}^d\left[\sum_{u=0}^\infty\beta^uh_i(u) 
			+ \sum_{u=0}^\infty\beta^u(1-\beta)\underline M_i(u)\right],$$
		which is integrable thanks to the boundedness of the sequences $h_i(\cdot)$ and 
		$\underline M_i(\cdot)$. Therefore, $\Lambda(\infty)$ is almost surely absolutely convergent, 
		and the process $\{\Lambda(t)\}_{t\in\N_0}$ is uniformly integrable.
	\end{proof}
	
	\begin{proof}[Proof of \eqref{itm:t1}]
		Each process $U_i(\cdot)$ of Proposition \ref{prop:U-mtg} is an $\FF_i$-martingale, so the process 
		$K(\cdot)$ of \eqref{eq:t1} is an $\FF(\tilde T^*)$-martingale, by Lemma \ref{lem:K-mtg}. In the
		notation of \eqref{eq:main-2}, we may write $K(\cdot)=\Lambda(\cdot)+W(\cdot)$, so it suffices to
		show that the families $\{\Lambda(t)\}_{t\in\N_0}$ and $\{W(t)\}_{t\in\N_0}$ are each uniformly 
		integrable; we have already shown that $\{\Lambda(t)\}_{t\in\N_0}$ is uniformly integrable 
		in the proof of \eqref{itm:t2}.
		
		Fix $t\in\N_0$ and observe, by an argument analogous to that of \eqref{eq:main-3} and \eqref{eq:main-4}, that
		$$W(t) = \sum_{u=\rho(k;\tilde T^*)}^{t-1}\beta^{u-T_j^*(u)}
			\Big(\beta^{T_j^*(u+1)}Y_j(u+1) - \beta^{T_j^*(u)}Y_j(u)\Big) = \beta^tY_j(t),$$
		where $k$ is the largest integer for which $t>\rho(k;\tilde T^*)$, and $j\equiv j(t)$ is the unique project engaged at time $t$.
		
		Thus, we wish to show the uniform integrability of $\{\beta^tY_{j(t)}(t)\}_{t\in\N_0}$. For 
		every $s\in\N_0$ and $m\ge0$, we have by \eqref{eq:ost-decr-rep-2}, \eqref{eq:U-1}, and Remark 
		\ref{rmk:enlarge} that
		\begin{align*}
			\beta^t[V_i(s;m)-m] 
			&\le (1-\beta)\,\E\left[\sum_{\theta=s}^\infty\beta^\theta
				\underline M_i(s,\theta)\bigmid\F^i(s)\right]
			\le \E\left[\sum_{\theta=0}^\infty\beta^\theta h_i(\theta+1)\bigmid\F^i(s)\right]
		\end{align*}
		almost surely for each $i=1,\dots,d$. Since $Y_{j(t)}(t)$ is $\F(\tilde T^*(t))$-measurable, the tower rule yields
		$$\beta^tY_{j(t)}(t) \le \E\left[\sum_{\theta=0}^\infty \beta^\theta 
			h_{j(t)}(\theta+1)\bigmid\F(\tilde T^*(t))\right], \qquad\text{a.s.}$$
		As a consequence of the uniform boundedness of $\{h_i(\theta)\}_{\theta\in\N_0}$ over every $i=1,\dots,d$, the family 
		$$\left\{\E\left[\sum_{\theta=0}^\infty\beta^\theta h_{j(t)}(\theta+1)
			\bigmid \F(\tilde T^*(t))\right]\right\}_{t\in\N_0}$$
		is uniformly integrable \cite[Thm.~4.6.1]{Durrett2019}, hence the same holds for $\{W(t)\}_{t\in\N_0}$.
	\end{proof}
	
	Before proving \eqref{itm:t3}, we require one more technical lemma, characterizing time changes with respect to an allocation strategy.
		
	\begin{lemma}\label{lem:time-change}
	Let $\tilde T\in\A(\tilde s)$ be an allocation strategy. Suppose $f_i:\N_0\to[0,\infty)$, $i=1,\dots,d$, are non-decreasing, and set $F(t):=\sum_{i=1}^df_i(T_i(t))$. Then for any $\theta<t$ and $\beta\in(0,1)$, we have
	\begin{equation}\label{eq:time-change}
	\sum_{i=1}^d\sum_{u=\theta}^{t-1} \beta^{F(u)-f_i(T_i(u))} \left(\beta^{f_i(T_i(u))} - \beta^{f_i(T_i(u+1))}\right) =  \beta^{F(\theta)} - \beta^{F(t)} \ge 0.
	\end{equation}
	\end{lemma}
	\begin{proof}
	For ease of notation, set $F_i(t)=f_i(T_i(t))$, and define the finite difference operator $\Delta f(t) := f(t+1)-f(t)$. The inner sum, via summation by parts, equals
	$$\beta^{F(t-1)+\Delta F_i(t-1)} - \beta^{F(\theta)} + \sum_{u=\theta}^{t-2} \beta^{F(u)}\left(\beta^{\Delta F_i(u)} - \beta^{\Delta F(u)}\right).$$
	Since allocation strategies can engage only one project at a time, then for each $u=\theta,\dots,t-1$, we must have either $\Delta F(u)=\Delta F_i(u)=0$ for all $i=1,\dots,d$, or there exists a unique $j\equiv j(u)$ for which $\Delta F(u) = \Delta F_j(u)>0$ while $\Delta F_i(u)=0$ for $i\neq j$. Thus, letting $\theta\le u_0<\dots<u_n<t-1$ be the points at which $\Delta F(u)>0$, we obtain
	\begin{multline*}
	\sum_{i=1}^d\sum_{u=\theta}^{t-2} 
		\beta^{F(u)}\left(\beta^{\Delta F_i(u)} - \beta^{\Delta F(u)}\right) \\
	= \sum_{u=\theta}^{t-2} \sum_{\substack{i=1\\i\neq j(u)}}^d \beta^{F(u)}\left(1-\beta^{\Delta F(u)}\right)
	= (d-1)\sum_{k=0}^n \beta^{F(u_k)}\left(1-\beta^{\Delta F(u_k)}\right).
	\end{multline*}
	Similarly, we find $\sum_{i=1}^d \beta^{F(t-1)+\Delta F_i(t-1)} = (d-1)\beta^{F(t-1)} + \beta^{F(t)}$. Consequently, the left-hand side of \eqref{eq:time-change} is given by
	$$\beta^{F(t)}-\beta^{F(\theta)} + (d-1)\left[\beta^{F(t-1)} - \beta^{F(\theta)} + \sum_{k=0}^n\left(\beta^{F(u_k)}-\beta^{F(u_{k+1})}\right)\right],$$
	using the fact that $F(u_k)+\Delta F(u_k) = F(u_k+1) = F(u_{k+1})$, with the convention that $F(u_{n+1}) = F(t-1)$. Now the summation telescopes and equals $\beta^{F(u_0)}-\beta^{F(u_{n+1})}$, but by definition $F(u_0)=F(\theta)$ and $F(u_{n+1})=F(t-1)$. Thus the term inside the square brackets is zero, completing the proof.
	\end{proof}
	
\begin{proof}[Proof of \eqref{itm:t3}]
	For each project $i=1,\dots,d$, let $\{V_i(t)\}_{t=s_i}^\infty$ and $\{Z_i(t)\}_{t=s_i}^\infty$ be the corresponding value processes and Snell envelopes, respectively, in the notation of \eqref{eq:osp} and \eqref{eq:snell} with $m=0$. Set
	\begin{equation}\label{eq:X1}
	X_1(t;\tilde T) := \sum_{i=1}^d\sum_{u=0}^{t-1} \beta^{u-T_i(u)}\Big(Z_i(T_i(u+1)) - Z_i(T_i(u))\Big)
	\end{equation}
	as well as
	\begin{equation}\label{eq:X2}
	X_2(t;\tilde T) := \beta^t\widehat\Phi(\tilde T(t))
	- \sum_{i=1}^d\sum_{u=0}^{t-1} \beta^{u-T_i(u)}\Big(\beta^{T_i(u+1)}V_i(T_i(u+1)) - \beta^{T_i(u)}V_i(T_i(u))\Big).
	\end{equation}
	It follows from \eqref{eq:Z-1} that $Q(t;\tilde T) = X_1(t;\tilde T)+X_2(t;\tilde T)$, so it suffices to show $\{X_1(t;\tilde T)\}_{t\in\N_0}$ and $\{X_2(t;\tilde T)\}_{t\in\N_0}$ are non-negative, bounded $\FF(\tilde T)$-supermartingales. Since Snell envelopes are non-negative $\FF_i$-supermartingales, the former follows from Lemma \ref{lem:K-mtg}. For the latter, observe that \eqref{eq:dap-decr-V} and \eqref{eq:tau} imply
	$$\beta^t\widehat\Phi(\tilde T(t)) = \E\left[\int_0^\infty\left(\beta^t-\beta^{\sum_{i=1}^d\sigma_i(T_i(t);m)}\right)\dd m\bigmid\F(\tilde T(t))\right]$$
	while \eqref{eq:ost-decr-rep-1} in conjunction with Proposition \ref{prop:enlarge} yields, for each $i=1,\dots,d$,
	$$\beta^{T_i(t)}V_i(T_i(t)) = \E\left[\int_0^\infty\left(\beta^{T_i(t)} - \beta^{\sigma_i(T_i(t);m)}\right)\dd m\bigmid\F^i(T_i(t))\right].$$
	We aim to show, for every $\theta<t$, that $\E[X_2(t;\tilde T)\mid\F(\tilde T(\theta))] \le X_2(\theta;\tilde T)$. This is equivalent, upon rearranging terms in \eqref{eq:X2}, to showing
	\begin{multline}\label{eq:t3-big}
	\E\left[\int_0^\infty \left(\beta^{\sum_{i=1}^d\sigma_i(T_i(t);m)} - \beta^{\sum_{i=1}^d\sigma_i(T_i(\theta);m)} + \phantom{\sum_{i=1}^d} \right.\right. \\
	\left.\left. + \sum_{i=1}^d\sum_{u=\theta}^{t-1} \beta^{u-T_i(u)} \left(\beta^{\sigma_i(T_i(u);m)} - \beta^{\sigma_i(T_i(u+1);m)}\right)\right)\dd m\bigmid\F(\tilde T(\theta))\right] \\
	\ge \E\left[\int_0^\infty\left(\beta^t - \beta^\theta + \sum_{i=1}^d\sum_{u=\theta}^{t-1} \beta^{u-T_i(u)}\left(\beta^{T_i(u)} - \beta^{T_i(u+1)}\right)\right)\dd m\bigmid\F(\tilde T(\theta))\right].
	\end{multline}
	The integrand on the right-hand side is identically zero by Lemma \ref{lem:time-change}, setting each $f_i$ to be the identity. On the other hand, observe that $u-T_i(u) = \sum_{j\neq i} T_j(u) \le \sum_{j\neq i}\sigma_j(T_j(u);m)$ almost surely, so the integrand on the left-hand side dominates
	\begin{multline*}
	\beta^{\sum_{i=1}^d\sigma_i(T_i(t);m)} - \beta^{\sum_{i=1}^d\sigma_i(T_i(\theta);m)}\ + \\
	+ \sum_{i=1}^d\sum_{u=\theta}^{t-1} \beta^{\sum_{j\neq i}\sigma_j(T_j(u);m)} \left(\beta^{\sigma_i(T_i(u);m)} - \beta^{\sigma_i(T_i(u+1);m)}\right)
	\end{multline*}
	which is non-negative, again by Lemma \ref{lem:time-change} with $f_i(t) := \sigma_i(t;m)$. This completes the proof of \eqref{itm:t3} and thus of Theorem \ref{thm:main} as well.
	\end{proof}

	%
	%
	
	\section{The Whittle reduction}\label{s:whittle}
	
	In this section, we show that the value of the dynamic allocation problem \eqref{eq:dap} has {\em yet another representation}, in addition to the two in \eqref{eq:dap-V}, known as the `Whittle reduction' after \cite{Whittle1980}. This representation was derived in the non-Markovian setting in \cite{ElKaroui1993}, which our proof essentially follows, except we demonstrate that one can dispense with the independence condition assumed in that work. The key point is that the equivalence of the representation \eqref{eq:whittle-V} presented below and the representations in \eqref{eq:dap-V} was previously established in a rather straightforward way, but relying on the independence of filtrations, whereas Theorem \ref{thm:main} constitutes a large detour required to prove \eqref{eq:dap-V}. Finally, we remark that Theorem \ref{thm:dap-indep} provides an alternative proof of the optimality of index-type strategies.
	
	We appeal to the following characterization of the random field $\Phi(\,\cdot\,;M)$ of \eqref{eq:dap}.
	\begin{lemma}\label{lem:bellman}
		For every given  $M\ge0$, the random field $\Phi(\,\cdot\,;M)$ of \eqref{eq:dap} is the 
		unique bounded solution of the {\bf Bellman Dynamic Programming Equation}
		\begin{equation}\label{eq:bellman}
			\Phi(\tilde s;M) = \max\Bigg[M, \max_{1\le j\le d}\Big(h_j(s_j+1)
				+ \beta\,\E[\Phi(\tilde s+\tilde e_j;M) \mid \F(\tilde s)]\Big)\Bigg], 
				\qquad\text{a.s.}
		\end{equation}
	\end{lemma}
	
	Analogous to the dynamic programming equation \eqref{eq:dp} for the problem of optimal stopping, the Bellman equation states that the maximum future expected discounted reward in the dynamic allocation problem can be attained by one of two ways: either by retiring immediately and collecting the retirement reward $M$; or by continuing for one additional unit of time, engaging a particular project, and evaluating the maximum future conditional expected reward at that point in time. The proof, which can be found in the \hyperref[appendix]{Appendix}, is from \cite{Mandelbaum1981, Mandelbaum1986}, and does not rely on the assumptions that the filtrations $\FF_1,\dots,\FF_d$ are independent.
		
	\vspace{2mm}
	Following \cite{Whittle1980}, we derive now an explicit expression for the random field $\Phi(\tilde s;M)$ of \eqref{eq:dap}. Recall that $K$ is a fixed upper bound on the maximum discounted reward that can possibly be accumulated from any single project, and that $M_i(t)$ denotes the Gittins index of the $i$th project at time $t$.
	
	\begin{theorem}\label{thm:dap-indep}
		Assume that the multi-parameter filtration $\FF$ satisfies the conditions \eqref{itm:F1}, \eqref{itm:F2}, and \eqref{itm:F4}. The random field
		\begin{equation}\label{eq:whittle-V}
			F(\tilde s;M) := \left.\begin{cases}
				K-\int_M^K\left(\prod_{i=1}^d\frac{\partial^+}{\partial m}V_i(s_i;m)\right) \dd m 
					&\text{if }\, 0\le M<K \\
				M\hfill &\text{if }\, M\ge K
			\end{cases}\right\},\qquad \tilde s\in\N_0^d
		\end{equation}
		satisfies the Bellman Dynamic Programming Equation \eqref{eq:bellman}. 
		
		In particular, for every 
		$M\ge0$, $\tilde s\in\N_0^d$, and $i=1,\dots,d$, we have
		\begin{equation}\label{eq:whittle-1}
			F(\tilde s;M)\ge M
		\end{equation}
		and 
		\begin{equation}\label{eq:whittle-2}
			F(\tilde s;M)\ge h_i(s_i+1)+\beta\,\E[F(\tilde s+\tilde e_i;M)\mid\F(\tilde s)].
		\end{equation}
		
		Moreover, recalling the notation of \eqref{eq:M-cal}, we have
		\begin{equation}\label{eq:whittle-3}
			F(\tilde s;M) = M \quad\text{on the event}\quad
			\left\{M\ge\M(\tilde s)\right\}
		\end{equation}
		and
		\begin{equation}\label{eq:whittle-4}
			\begin{split}
			F(\tilde s;M) = h_i(s_i+1)+\beta\,\E[F(\tilde s+\tilde e_i;M)\mid\F(\tilde s)] \\
				\text{on the event}\quad\left\{M<\M(\tilde s)=:M_i(s_i)\right\}.
			\end{split}
		\end{equation}
	\end{theorem}
	\begin{proof}
		By Lemma \ref{lem:V-rdiff-props}, $0\le\frac{\partial^+}{\partial m}V_i(s_i;m)\le1$ holds for each
		$i=1,\dots,d$, hence $F(\tilde s;M)\ge M$ when $M<K$, which proves \eqref{eq:whittle-1}. Also by Lemma
		\ref{lem:V-rdiff-props}, $\frac{\partial^+}{\partial m}V_i(s_i;m)=1$ holds if $m\ge M_i(s_i)$, 
		hence $\prod_{i=1}^d\frac{\partial^+}{\partial m}V_i(s_i;m)=1$ if
		$m\ge\M(\tilde s)$. Then the integrand in \eqref{eq:whittle-V} is identically 
		equal to one when $M\ge\M(\tilde s)$, which proves \eqref{eq:whittle-3}.
		
		To prove \eqref{eq:whittle-2}, first consider the case where $M\ge K$. From the bound 
		$0\le h_i\le K(1-\beta)$, we easily obtain
		$$h_i(s_i+1)+\beta\,\E[F(\tilde s+\tilde e_i;M)\mid\F(\tilde s)] \le K(1-\beta)+\beta M \le M 
			= F(\tilde s;M).$$
			
		Now suppose $0\le M<K$. Fix $i\in\{1,\dots,d\}$ and denote the vector given by removing the 
		$i$th coordinate of $\tilde s$ by $\tilde s^{(i)}:=(s_1,\dots,s_{i-1},s_{i+1},\dots,s_d)
		\in\N_0^{d-1}$. Consider the mapping
		$$m\mapsto P_i\left(\tilde s^{(i)};m\right) 
			:= \prod_{\substack{j=1\\j\neq i}}^d \frac{\partial^+}{\partial m} V_j(s_j;m), 
			\qquad m\ge0.$$
		Observe that, by Lemma \ref{lem:V-rdiff-props} and \eqref{eq:gittins-upper}, we have 
		$P_i(\tilde s^{(i)};K)=1$. Additionally, \eqref{eq:dp}, \eqref{eq:sigma-equivs}, and 
		\eqref{eq:gittins-upper} imply that $V_i(s_i;K)=K$. Integrating by parts then yields
		\begin{align}\label{eq:whittle-5}
			F(\tilde s;M)
			&= K - \int_M^K P_i\left(\tilde s^{(i)};m\right) 
				\cdot \frac{\partial^+}{\partial m}V_i(s_i;m)\,\dd m \nonumber \\
			&= V_i(s_i;M)P_i\left(\tilde s^{(i)};M\right)
				+ \int_M^K V_i(s_i;m)\,\dd_mP_i\left(\tilde s^{(i)};m\right).
		\end{align}
		Therefore,
		$$F(\tilde s+\tilde e_i;M) = V_i(s_i+1;M)P_i\left(\tilde s^{(i)};M\right)
			+ \int_M^K V_i(s_i+1;m)\,\dd_mP_i\left(\tilde s^{(i)};m\right).$$
		Since $V_i(s_i+1;,)\in\F_i(s_i+1)$, the \eqref{itm:F4} condition in the form \eqref{eq:F4-1} yields $\E[V_i(s_i+1;m)\mid
		\F(\tilde s)] = \E[V_i(s_i+1;m)\mid\F_i(s_i)]$ for all $m\ge0$. Since
		$P_i(\tilde s^{(i)};M)$ is $\F(\tilde s)$-measurable, taking conditional expectations and 
		applying Fubini's theorem yields
		\begin{multline}\label{eq:whittle-6}
			\E[F(\tilde s+\tilde e_i;M)\mid\F(\tilde s)] = P_i\left(\tilde s^{(i)};M\right)\E[V_i(s_i+1;M)\mid\F_i(s_i)]\ + \\
			+ \int_M^K \E[V_i(s_i+1;m)\mid\F_i(s_i)]\,\dd_m 
				P_i\left(\tilde s^{(i)};m\right).
		\end{multline}
		Now define 
		$$\phi_i(t;m) := V_i(t;m) - h_i(t+1) - \beta\,\E[V_i(t+1;m)\mid\F_i(t)].$$
		By \eqref{eq:whittle-5}, \eqref{eq:whittle-6}, and the fact that $P_i(\tilde s^{(i)};K)=1$, we obtain the identity
		\begin{equation}\label{eq:whittle-7}
		\begin{split}
			\phi_i(s_i;M)P_i &\left(\tilde s^{(i)}; M\right) 
				+ \int_M^K \phi_i(s_i;m)\,\dd_mP_i\left(\tilde s^{(i)};m\right) \\
			&= F(\tilde s;M)-h_i(s_i+1)-\beta\,\E[F(\tilde s+\tilde e_i;M)\mid\F(\tilde s)].
		\end{split}
		\end{equation}
		The left-hand side of \eqref{eq:whittle-7} is non-negative by \eqref{eq:dp}, so the
		right-hand side is non-negative as well, proving \eqref{eq:whittle-2}.
		
		It remains to show \eqref{eq:whittle-4}. On $\{M<M_i(s_i):=\M(\tilde s)\}$, we have $M<K$
		by \eqref{eq:gittins-upper}, so it suffices to show that the left-hand side of \eqref{eq:whittle-7} is zero in 
		this event.
		
		Let $m\in[M,K]$ be given. On $\{m<M_i(s_i)\}$, we have by \eqref{eq:sigma-equivs} that $V_i(s_i;m)>m$,
		hence \eqref{eq:dp} implies $\phi_i(s_i;m)=0$. In particular, $\phi_i(s_i;M)=0$.
		On the other hand, we have by definition that $m\ge M_j(s_j)$ holds for each $j=1,\dots,d$ on 
		$\{m\ge M_i(s_i)\}$. Lemma \ref{lem:V-rdiff-props} then implies $P_i(\tilde s^{(i)};m)=1$, hence
		$$\int_M^K\phi_i(s_i;m)\,\dd_mP_i\left(\tilde s^{(i)};m\right)
			= \int_M^{M_i(s_i)}0\,\dd_mP_i\left(\tilde s^{(i)};m\right)
				+ \int_{M_i(s_i)}^K \phi_i(s_i;m)\,\dd_m(1) = 0,$$
		which proves \eqref{eq:whittle-4}.
	\end{proof}

\appendix

\section{Proofs}\label{appendix}

	\subsection{Proofs of results in Section \ref{s:osp}}
	\begin{proof}[Proof of Lemma \ref{lem:sigma-props}]
		Let $\phi(t;m) := \beta^t(V(t;m) - m)$. Each of the following claims holds almost surely. We have from \eqref{eq:V-lower} that $\phi\ge0$, and by \eqref{eq:osp}, we can express $\phi$ via a telescoping sum as
		\begin{align}\label{phi decrease}
			\phi(t;m) 
			&= \esssup_{\tau\in\s(t)}\,\E\left[\sum_{u=t}^{\tau-1} \beta^u (h(u+1) - m(1-\beta))
				\bigmid \F(t) \right].
		\end{align}
		It is then clear that $m\mapsto\phi(t;m)$ is strictly decreasing. Let $m_1>m_2$ be given and note 
		that $V(\sigma(t;m_2);m_2) = m_2$ by \eqref{eq:sigma-V}, hence $\phi(\sigma(t;m_2);m_2) = 0$. By
		monotonicity we have
		\begin{equation}\label{phi sigma = 0}
			0 \le \phi(\sigma(t;m_2);m_1) \le \phi(\sigma(t;m_2);m_2) = 0
		\end{equation}
		which implies $V(\sigma(t;m_2);m_1) = m_1$. By definition $\sigma(t;m_1)$ is the 
		smallest integer $\theta$ for which $V(\theta;m_1) = m_1$, so $\sigma(t;m_1) \le \sigma(t;m_2)$.
		Thus $m\mapsto\sigma(t;m)$ is decreasing.
		
		In showing right-continuity, we have from the just-established decrease that, for fixed $m\ge0$, 
		$$\sigma_* := \lim_{\delta\downarrow0} \sigma(t;m+\delta) \le \sigma(t;m)$$
		where the limit exists because $\sigma$ is bounded from below by $t$.
		Fix $\delta>0$ and observe that for integer $k>1/\delta$, we have 
		$\phi(\sigma(t;m+1/k);m+\delta) = 0$ by analogy with \eqref{phi sigma = 0}. 
		Since $\sigma(t;\,\cdot\,)$ takes integer
		values $\sigma(t;m+1/k) = \sigma_*$ holds for all large enough $k$, hence $\phi(\sigma_*;
		m+\delta)=0$ for all $\delta>0$. Now we see from \eqref{eq:osp} that
		$m\mapsto V(t;m)$ is continuous, so $m\mapsto\phi(t;m)$ is also continuous. 
		Taking $\delta\downarrow0$, we get $\phi(\sigma_*;m)=0$. We conclude by the minimality of 
		$\sigma(t;m)$ that $\sigma(t;m)\le\sigma_*$, and this completes the proof.
	\end{proof}
	
	\begin{proof}[Proof of Lemma \ref{lem:V-rdiff}]
		Increase and convexity follow immediately from the definition \eqref{eq:osp}. 
		Again the following claims hold almost surely.
		
		For \eqref{eq:V-rdiff}, fix $m\ge0$. From the definition and Lemma \ref{lem:sigma-props}, 
		we have for any $\delta>0$ the inequalities 
		$$t\le\sigma(t;m+\delta)\le\sigma(t;m).$$
		From Proposition \ref{prop:Z-mtg} and \eqref{eq:Z-1}, we know the stopped process
		$Z(\,\cdot\,\wedge\sigma(t;m);m)$ is a non-negative martingale, hence we may apply 
		the optional sampling theorem \cite[Thm.~4.8.4]{Durrett2019} to obtain
		$$Z(t;m) = \E[Z(\sigma(t;m+\delta);m)\mid\F(t)].$$
		Substituting, and using \eqref{eq:Z-1}, this equality becomes
		\begin{multline*}
			\sum_{u=0}^{t-1} \beta^uh(u+1) + \beta^tV(t;m) \\
			= \E\left[\sum_{u=0}^{\sigma(t;m+\delta)-1}\beta^uh(u+1) + \beta^{\sigma(t;m+\delta)}
				V(\sigma(t;m+\delta);m) \bigmid \F(t)\right].
		\end{multline*}
		Isolating $V(t;m)$ and using the inequality \eqref{eq:V-lower}, we obtain
		\begin{align}\label{l2 step}
			V(t;m) 
			&\ge \E\left[\sum_{u=t}^{\sigma(t;m+\delta)-1}\beta^{u-t}h(u+1) +
				m\beta^{\sigma(t;m+\delta)-t} \bigmid \F(t)\right].
		\end{align}
		We recall at this stage \eqref{eq:sigma-optimal}, which, upon replacing $m$ with $m+\delta$, 
		yields
		$$V(t;m+\delta) = \E\left[\sum_{u=t}^{\sigma(t;m+\delta)-1}\beta^{u-t}h(u+1) 
			+ (m+\delta)\beta^{\sigma(t;m+\delta)-t}\mid\F(t)\right].$$
		From this in conjunction with \eqref{l2 step}, we have
		\begin{align}\label{l2 limsup}
			\limsup_{\delta\downarrow0}\frac{V(t;m+\delta)-V(t;m)}{\delta}
			\le \limsup_{\delta\downarrow0}\,\E[\beta^{\sigma(t;m+\delta)-t}\mid\F(t)]
			\le \E[\beta^{\sigma(t;m)-t}\mid\F(t)],
		\end{align}
		where the second inequality holds on the strength of Fatou's lemma and Lemma \ref{lem:sigma-props}.
		
		On the other hand, we have 
		$$Z(t;m+\delta)\ge \E[Z(\sigma(t;m);m+\delta)\mid\F(t)]$$
		by applying the optional sampling theorem once again on the supermartingale 
		$Z(\,\cdot\,;m+\delta)$. By analogy with \eqref{l2 step}, we obtain
		\begin{align*}
			V(t;m+\delta) 
			&\ge \E\left[\sum_{u=t}^{\sigma(t;m)-1}\beta^{u-t}h(u+1) + (m+\delta)\beta^{\sigma(t;m)-t}
				\bigmid \F(t)\right],
		\end{align*}
		From this, in conjunction with \eqref{eq:sigma-optimal}, it follows that
		\begin{equation}\label{l2 liminf}
			\liminf_{\delta\downarrow0} \frac{V(t;m+\delta)-V(t;m)}{\delta} \ge 
			\E[\beta^{\sigma(t;m)-t}\mid\F(t)].
		\end{equation}
		The inequalities \eqref{l2 limsup} and \eqref{l2 liminf} thus prove \eqref{eq:V-rdiff}.
	\end{proof}

	\subsection{Proofs of results in Section \ref{s:gittins}}
	
	\begin{proof}[Proof of \eqref{eq:gittins-forwards}]
		Since $\sum_{u=0}^\infty\beta^u=(1-\beta)^{-1}$ then we can write the left-hand side of
		\eqref{eq:gittins-forwards} in terms of \eqref{eq:osp}, i.e.,
		\begin{equation}\label{fair charge}
			(1-\beta)M(t)
			= \esssup\left\{X\in\M(t)\colon \esssup_{\tau\in\s(t)}\,\E\left[\sum_{u=t}^{\tau-1}\beta^{u-t}
				\Big(h(u+1)-X\Big)\bigmid\F(t)\right]\ge0\right\}.
		\end{equation}
		For every $\tau\in\s(t+1)$, let
		\begin{equation}\label{fair charge'}
			\lambda(\tau) := \frac{\E\left[\sum_{u=t}^{\tau-1}\beta^uh(u+1)\bigmid\F(t)\right]}
				{\E\left[\sum_{u=t}^{\tau-1}\beta^u\bigmid\F(t)\right]},
			\tag{\ref*{fair charge}$'$}
		\end{equation}
		noting that $\lambda(\tau)\in\M(t)$ since $h$ is positive. Additionally, let 
		$$\lambda^*:=\esssup_{\tau\in\s(t+1)}\lambda(\tau)\in\M(t),$$
		so there exists a sequence of stopping times
		$\tau_n\in\s(t+1)$, $n\in\N$, such that $\lambda(\tau_n)\to\lambda^*$ almost surely. Notice 
		from \eqref{fair charge'} that
		\begin{equation}\label{fair charge''}
			\E\left[\sum_{u=t}^{\tau-1} \beta^{u-t}
				\Big(h(u+1)-\lambda(\tau)\Big)\bigmid\F(t)\right] = 0, 
				\quad\tau\in\s(t+1),\qquad \text{a.s.},
			\tag{\ref*{fair charge}$''$}
		\end{equation}
		which yields
		\begin{equation}\label{FI1}
			\E\left[\sum_{u=t}^{\tau-1}\beta^{u-t}\Big(h(u+1)-\lambda^*\Big)\bigmid\F(t)\right] \le 0,
			\quad\tau\in\s(t+1),\qquad \text{a.s.},
		\end{equation}
		and, by the bounded convergence theorem,
		\begin{equation}\label{FI2}
			\lim_{n\to\infty}
				\E\left[\sum_{u=t}^{\tau_n-1}\beta^{u-t}\Big(h(u+1)-\lambda^*\Big)\bigmid\F(t)\right]
				= 0, \qquad\text{a.s.}
		\end{equation}
		
		On the other hand, \eqref{FI1} implies 
		\begin{equation}\label{FI3}
			\esssup_{\tau\in\s(t)}\,\E\left[\sum_{u=t}^{\tau-1}\beta^{u-t}\Big(h(u+1)-\lambda^*\Big)
				\bigmid\F(t)\right] \le 0, \qquad\text{a.s.},
		\end{equation}
		hence $\lambda^*$ is almost surely an upper bound of the set
		$$\left\{X\in\M(t)\colon \esssup_{\tau\in\s(t)}\,\E\left[\sum_{u=t}^{\tau-1}\beta^{u-t}
			\Big(h(u+1)-X\Big)\bigmid\F(t)\right]\ge0\right\}.$$
		But $\lambda^*$ is an element of this set, modulo events of null measure, as follows from 
		\eqref{FI2} and \eqref{FI3}. It follows from \eqref{fair charge} that $\lambda^*=(1-\beta)M(t)$ 
		almost surely, as desired.
	\end{proof}
	
	\subsection{Proofs of results in Section \ref{s:strategies}}
	
	\begin{proof}[Proof of Proposition \ref{prop:excursion}]
	Without loss of generality, assume $\tilde s = \tilde 0$. Recall that \eqref{eq:sync} is equivalent to the lower envelope index-type property \eqref{eq:lower-index-type}. 
	
	($\Longrightarrow$) We proceed by induction on $t$.  At $t=0$, we have $M_j(T_j(0)) = \underline M_j(T_j(0))$ for every $j=1,\dots,d$, hence \eqref{eq:lower-index-type} implies \eqref{eq:index-type} at $t=0$.
	
	Assume then that \eqref{eq:index-type} is satisfied for some $t\ge0$, and that $M_j(T_j(t))=\underline M_j(T_j(t))$ for all $j\neq i$ on the event $\{\tilde T(t+1)=\tilde T(t)+\tilde e_i\}$. Now suppose $\tilde T(t+1)=\tilde T(t)+\tilde e_i$ and $\tilde T(t+2)=\tilde T(t+1)+\tilde e_k$ for some $i,k$.
	
	{\it Case 1}: $M_i(T_i(t+1)) = \underline M_i(T_i(t+1))$. Then the same equality holds for every $j=1,\dots,d$ by the inductive hypothesis, so \eqref{eq:lower-index-type} implies $M_k(T_k(t))=\M(\tilde s)$, in the notation of \eqref{eq:M-cal}, as desired.
	
	{\it Case 2}: $M_i(T_i(t+1)) > \underline M_i(T_i(t+1)) = \underline M_i(T_i(t))$. Then \eqref{eq:excursion} implies $i=k$, and \eqref{eq:lower-index-type} in conjunction with the inductive hypothesis yields
	$$M_i(T_i(t+1)) > \underline M_i(T_i(t)) = \underline\M(\tilde T(t)) \ge \underline M_j(T_j(t))
		= M_j(T_j(t)) = M_j(T_j(t+1))$$
	for every $j\neq i$. This completes the induction.
	
	($\Longleftarrow$) We follow a similar strategy. Assume that \eqref{eq:sync} and \eqref{eq:excursion} hold for some $t\ge0$, and that $M_j(T_j(t)) = \underline M_j(T_j(t))$ for all $j\neq i$ on the event $\{\tilde T(t+1)=\tilde T(t)+\tilde e_i\}$.
	
	{\it Case 1}: $M_i(T_i(t+1)) = \underline M_i(T_i(t+1))$. Clearly \eqref{eq:excursion} holds trivially for $t+1$. By the inductive hypothesis, $M_j(T_j(t+1))=\underline M_j(T_j(t+1))$ for every $j=1,\dots,d$. If $\tilde T(t+2)=\tilde T(t+1)+\tilde e_k$ for some $k$, then \eqref{eq:index-type} implies 
	$$\underline M_k(T_k(t+1)) = M_k(T_k(t+1)) = \M(\tilde T(t+1)) = \underline\M(\tilde T(t+1)),$$
	as desired.
	
	{\it Case 2}: $M_i(T_i(t+1)) > \underline M_i(T_i(t+1)) = \underline M_i(T_i(t))$. If $\tilde T(t+2)=\tilde T(t+1)+\tilde e_k$ for some $k\neq i$, then the inductive hypothesis yields
	$$M_k(T_k(t+1)) = \underline M_k(T_k(t)) \le \underline M_i(T_i(t)) = \underline M_i(T_i(t+1)) < M_i(T_i(t+1)),$$
	contradicting \eqref{eq:index-type}. Thus $k=i$, which shows \eqref{eq:excursion} for $t+1$. In particular, we have 
	$$\underline M_k(T_k(t+1)) = \underline M_i(T_i(t)) = \underline\M(\tilde T(t)) = \underline\M(\tilde T(t+1)),$$
	which is \eqref{eq:lower-index-type} for $t+1$.
	\end{proof}
	
	\begin{proof}[Proof of Proposition \ref{prop:rho-stopping}]
	We consider the latter claim and prove the former claim on the way. We proceed by induction on $k$; the base case is trivial. Suppose $\rho(k;\tilde T)$ is an $\FF(\tilde T)$-stopping time for some $k\in\N_0$. It suffices to stay within the event $\{t>\rho(k;\tilde T)\}$. In this case, we have
	\begin{equation}\label{eq:rho-1}
	\{\rho(k+1;\tilde T) = t\} 
		= \bigcup_{\tilde k\in\N_0^d} \left[\{\tilde T(t) = \tilde\varepsilon(\tilde k;\tilde s)\} \cap
		\bigcup_{u=0}^{t-1} \{\rho(k;\tilde T)=u,\,\tilde T(u)\neq\tilde\varepsilon(\tilde k;\tilde s)\}
		\right].
	\end{equation}
	Fixing $\tilde k\in\N_0^d$, we deal with each event in the above expression separately.
	
	Fix $\tilde r\in\N_0^d$. Since $\{\varepsilon_i(k_i;s_i)=r_i\}\in\F_i(r_i)\subset\F^i(r_i)$, then
	$$\{\tilde\varepsilon(\tilde k;\tilde s) = \tilde r\} = \bigcap_{i=1}^d \{\varepsilon_i(k_i;s_i)=r_i\}
		\in \bigcap_{i=1}^d \F^i(r_i) = \F(\tilde r),$$
	where the last equality follows by definition \eqref{eq:big-filtr}. Thus $\tilde\varepsilon(\tilde k;\tilde s)$ is an $\FF$-stopping point. Consequently,
	$$\{\tilde T(t)=\tilde r\}\cap\{\tilde T(t)=\tilde\varepsilon(\tilde k;\tilde s)\} = \{\tilde\varepsilon(\tilde k;\tilde s)=\tilde r\}\in\F(\tilde r), \qquad \tilde r\in\N_0^d,$$
	hence $\{\tilde T(t)=\tilde\varepsilon(\tilde k;\tilde s)\}\in\F(\tilde T(t))$ by definition.
	
	The other two events that remain are $\{\rho(k;\tilde T)=u\}$ and $\{\tilde T(u)\neq\varepsilon(\tilde k;\tilde s)\}$, for $u=0,\dots,t-1$. The latter event is $\F(\tilde T(t))$-measurable by closure under complements, while the former is $\F(\tilde T(u))$-measurable by the inductive hypothesis. But $\FF(\tilde T)$ is a filtration by Proposition \ref{prop:multiparam-filtr}, so $u<t$ implies our event is also $\F(\tilde T(t))$-measurable.
	
	We conclude that \eqref{eq:rho-1} is a countable union of $\F(\tilde T(t))$-measurable events, hence $\{\rho(k+1;\tilde T)=t\}\in\F(\tilde T(t))$. Since this holds for every $t\in\N_0$, $\rho(k+1;\tilde T)$ is an $\FF(\tilde T)$-stopping time.
	\end{proof}
	
	
	\subsection{Proofs of results in Section \ref{s:multiparam}}
	
	\newcommand{\h}{\mathcal H}
	\newcommand{\LL}{\mathcal L}
	\begin{proof}[Proof of Example \ref{ex:F4-indep}]
		We fix $\tilde s,\tilde r\in\N_0^d$ and prove that
		\begin{equation}\label{indep F4 1}
			\P(A\cap B\mid\F(\tilde s\wedge\tilde r)) 
				= \P(A\mid\F(\tilde s\wedge\tilde r))\,\P(B\mid\F(\tilde s\wedge\tilde r))
		\end{equation}
		holds for every $A\in\F(\tilde s)$ and $B\in\F(\tilde r)$.
		
		Recall that $\F(\tilde s)=\bigvee_i\F_i(s_i)$ and $\F(\tilde r)=\bigvee_i\F_i(r_i)$. It 
		suffices to show, by the monotone class theorem, that \eqref{indep F4 1} holds for 
		$A=\bigcap_iA_i$ and $B=\bigcap_iB_i$, for some $A_i\in\F_i(s_i)$ and $B_i\in
		\F_i(r_i)$. Since the $\sigma$-algebras $\F_i(s_i\vee r_i)$, $i=1,\dots,d$, are independent, 
		and $A_i\cap B_i\in\F_i(s_i\vee r_i)$, then
		$$\P(A\cap B\mid\F(\tilde s\wedge\tilde r)) 
			= \P\left(\bigcap_{i=1}^d(A_i\cap B_i)\bigmid\F(\tilde s\wedge\tilde r)\right)
			= \prod_{i=1}^d \E[\1_{A_i}\1_{B_i}\mid\F(\tilde s\wedge\tilde r)].$$
		Notice that at least one of $\1_{A_i}$ or $\1_{B_i}$ is 
		$\F(\tilde s\wedge\tilde r)$-measurable for each $i=1,\dots,d$, yielding
		$$\E[\1_{A_i}\1_{B_i}\mid\F(\tilde s\wedge\tilde r)] 
			= \E[\1_{A_i}\mid\F(\tilde s\wedge\tilde r)]\,\E[\1_{B_i}\mid\F(\tilde s\wedge\tilde r)],$$
		and we obtain again by independence
		\begin{align*}
			\prod_{i=1}^d\E[\1_{A_i}\1_{B_i}\mid\F(\tilde s\wedge\tilde r)] 
			&= \E\left[\prod_{i=1}^d\1_{A_i}\bigmid\F(\tilde s\wedge\tilde r)\right]
				\E\left[\prod_{i=1}^d\1_{B_i}\bigmid\F(\tilde s\wedge\tilde r)\right] \\
			&= \P(A\mid\F(\tilde s\wedge\tilde r))\,\P(B\mid\F(\tilde s\wedge\tilde r))
		\end{align*}
		as desired.
	\end{proof}
	
	\begin{proof}[Proof of Example \ref{ex:F4-sheet}]
		We show that the filtration \eqref{eq:sheet-filtr} satisfies \eqref{eq:F4-1}.
		Let  $R(\tilde s)$ denote the rectangle $\{1,\dots,s_1\}\times\dots\times\{1,\dots,s_d\}$ for 
		$\tilde s\in\N_0^d$. Fix $\tilde s,\tilde r\in\N_0^d$. We may write $R(\tilde s) = 
		R(\tilde s\wedge\tilde r)\cup U$ and $R(\tilde r)=R(\tilde s\wedge\tilde r)\cup V$ for some 
		$U,V\subset\N_0^d$ such that $U$, $V$, and $R(\tilde s\wedge\tilde r)$ are pairwise disjoint. 
		Then 
		\begin{equation}\label{Bsheet filtration 1}
			\F(\tilde s) = \F(\tilde s\wedge\tilde r)\vee\G(U) 
				= \sigma(F\cap G\colon F\in\F(\tilde s\wedge\tilde r),\,G\in\G(U))
		\end{equation}
		and
		\begin{equation}\label{Bsheet filtration 2}
			\F(\tilde r) = \F(\tilde s\wedge\tilde r)\vee\G(V)
				= \sigma(F\cap G\colon F\in\F(\tilde s\wedge\tilde r),\,G\in\G(V))
		\end{equation}
		where $\G(U):=\sigma(W(A)\colon A\subseteq U))$. By the construction of the white noise process, 
		$\F(\tilde s\wedge\tilde r)$, $\G(U)$, and $\G(V)$ are independent.
		
		Let $\xi$ be a bounded $\F(\tilde r)$-measurable function. Since 
		$\E[\xi\mid\F(\tilde s\wedge\tilde r)]$ is $\F(\tilde s)$-measurable, it suffices to show,
		again by the monotone class theorem, that
		\begin{equation}\label{Bsheet filtration 3}
			\E[\,\E[\xi\mid\F(\tilde s\wedge\tilde r)]\1_A]=\E[\xi\1_A]
		\end{equation}
		holds for every $A=F\cap G$, for some $F\in\F(\tilde s\wedge\tilde r)$ and 
		$G\in\G(U)$. In this case, $\xi$ and $\1_G$ are independent, and we obtain
		$$\E[\,\E[\xi\mid\F(\tilde s\wedge\tilde r)]\1_A] 
			= \E[\,\E[\xi\1_F\mid\F(\tilde s\wedge\tilde r)]]E[\1_G]
			= \E[\xi\1_F]\,\E[\1_G] 
			= \E[\xi\1_A]$$
		as desired.
	\end{proof}
	
	\subsection{Proofs of results in Section \ref{s:whittle}}
	
	\begin{proof}[Proof of Lemma \ref{lem:bellman}]
		We first show that $\Phi(\cdot;M)$, as defined in \eqref{eq:dap}, solves \eqref{eq:bellman}.
		Fix $\tilde s\in\N_0^d$, $M\ge0$, and $j=1,\dots,d$. By considering a policy where the 
		stopping time is identically zero, we have $\Phi(\tilde s;M)\ge M$.
		
		On the other hand, there exists a sequence of policies 
		$\Pi_n\in\p(\tilde s+\tilde e_j)$ such that 
		$$\Phi(\tilde s+\tilde e_j;M)=\lim_{n\to\infty}\E[\reward(\Pi_n;M)\mid\F(\tilde s+\tilde e_j)]$$
		almost surely. Since $h_i(\cdot)$ is bounded in $[0,K(1-\beta)]$ then
		$$\E[\Phi(\tilde s+\tilde e_j;M)\mid\F(\tilde s)] 
			= \lim_{n\to\infty} \E[\reward(\Pi_n;M)\mid\F(\tilde s)].$$
		Since $h_j(\cdot)$ is $\FF_j$-predictable then we have almost surely that
		\begin{multline*}
			h_j(s_j+1) + \beta\,\E[\Phi(\tilde s+\tilde e_j;M)\mid\F(\tilde s)] \\
			= \lim_{n\to\infty} \E[h_j(s_j+1)+\beta\reward(\Pi_n;M)\mid\F(\tilde s)]
				\le \Phi(\tilde s;M).
		\end{multline*}
		Here, we use the fact that the $\Pi_n$'s are policies for $\tilde s+\tilde e_j$, hence the integrands 
		$h_j(s_j+1)+\beta\reward(\Pi_n;M)$ are in fact the total discounted rewards of policies in 
		$\p(\tilde s)$. The inequality then follows from the maximality of $\Phi(\cdot)$. Because 
		this holds for arbitrary $j=1,\dots,d$, we obtain the inequality
		$$\Phi(\tilde s;M)\ge\max\left[M,\max_{1\le j\le d}\Big(h_j(s_j+1)
			+ \beta\,\E[\Phi(\tilde s+\tilde e_j;M)\mid\F(\tilde s)]\Big)\right], \qquad\text{a.s.}$$
			
		To obtain the reverse inequality, we write the total discounted reward for any given policy 
		$\Pi=(\tilde T,\tau)\in\p(\tilde s)$ as
		\begin{equation}\label{Bellman l1}
			\reward(\tilde T,\tau;M) = \reward(\tilde T,0;M)\1\{\tau=0\} 
			+ \reward(\tilde T,\tau\vee1;M)\1\{\tau\ge1\}.
		\end{equation}
		As before, we have $\reward(\tilde T,0;M)=M$. Also notice that on $\{\tilde T(1)=\tilde s+
		\tilde e_j,\,\tau\ge1\}$, we have
		\begin{multline*}
			\reward(\tilde T,\tau\vee1; M)
			= h_j(s_j+1)\ + \\
			+ \beta\left(\sum_{i=1}^d\sum_{t=0}^{(\tau\vee1)-2} 
				\beta^t h_i(T_i(t+2))\Big(T_i(t+2)-T_i(t+1)\Big)
				+ M\beta^{(\tau\vee1)-1}\right).
		\end{multline*}
		Let $\tilde T^j$ denote any allocation strategy in $\A(\tilde s+\tilde e_j)$ such that 
		$\tilde T^j(t) = \tilde T(t+1)$ holds on the event $\{\tilde T(1)=\tilde s+\tilde e_j\}$. Then
		\begin{multline*}
			\reward(\tilde T,\tau\vee1;M)\1\{\tilde T(1)=\tilde s+\tilde e_j,\,\tau\ge1\} \\
			= \left(h_j(s_j+1)+\beta\reward(\tilde T^j,(\tau\vee1)-1;M)\right)
				\1\{\tilde T(1)=\tilde s+\tilde e_j,\,\tau\ge1\}.
		\end{multline*}
		We may thus write \eqref{Bellman l1} as
		\begin{multline*}
			\reward(\Pi;M) \\
			= M\1\{\tau=0\} + \sum_{j=1}^d\left(h_j(s_j+1)+\beta\reward(\tilde T^j,(\tau\vee1)-1;M)\right)
				\1\{\tilde T(1)=\tilde s+\tilde e_j,\,\tau\ge1\}.
		\end{multline*}
		It is clear that $(\tilde T^j,(\tau\vee1)-1)$ is a policy for $\tilde s+\tilde e_j$. By the 
		maximality of $\Phi$, we obtain
		\begin{align*}
			\E[\reward(\Pi;M)\mid\F(\tilde s)] 
			&\le \max\left[M, \max_{1\le j\le d}\Big(h_j(s_j+1)
				+ \beta\,\E[\Phi(\tilde s+\tilde e_j;M)\mid\F(\tilde s)]\Big)\right].
		\end{align*}
		This holds for an arbitrary policy $\Pi\in\p(\tilde s)$, hence
		$$\Phi(\tilde s;M) \le \max\left[M, \max_{1\le j\le d}\Big(h_j(s_j+1)
			+ \beta\,\E[\Phi(\tilde s+\tilde e_j;M)\mid\F(\tilde s)]\Big)\right]$$
		as desired.
		
		It remains to show uniqueness. Suppose $\Psi(\tilde s;M)$ is another bounded solution to 
		\eqref{eq:bellman}. Let $X(\tilde s;M)=|\Phi(\tilde s;M)-\Psi(\tilde s;M)|$.
		
		Suppose $\Phi(\tilde s;M)>M$. Then there exists some $k=1,\dots,d$ such that 
		$$\Phi(\tilde s;M) = h_k(s_k+1)+\beta\,\E[\Phi(\tilde s+\tilde e_k;M)\mid\F(\tilde s)].$$
		By \eqref{eq:bellman} applied to $\Psi$, we have
		\begin{align*}
			\Phi(\tilde s;M)-\Psi(\tilde s;M) 
			&\le \beta\,\E[\Phi(\tilde s+\tilde e_k;M)-\Psi(\tilde s+\tilde e_k;M)\mid\F(\tilde s)] \\
			&\le \max_{1\le j\le d}\beta\,\E[X(\tilde s+\tilde e_j;M)\mid\F(\tilde s)].
		\end{align*}
		If instead $\Phi(\tilde s;M)=M$ then \eqref{eq:bellman} again yields
		$$\Phi(\tilde s;M)-\Psi(\tilde s;M)\le 0
			\le \max_{1\le j\le d}\beta\,\E[X(\tilde s+\tilde e_j;M)\mid\F(\tilde s)].$$
		The same inequality holds upon interchanging $\Phi$ and $\Psi$, hence 
		$$X(\tilde s;M)\le\max_{1\le j\le d}\beta\,\E[X(\tilde s+\tilde e_j;M)\mid\F(\tilde s)].$$
		By iterating this inequality and applying the tower rule, we obtain
		$$0\le X(\tilde s;M) 
			\le \beta^n\max_{1\le j_1,\dots,j_n\le d}\E[X(\tilde s+e_{j_1}+\dots+e_{j_n};M)
			\mid\F(\tilde s)]$$
		for every $n\in\N$. Since $X$ is bounded due to the boundedness of $\Phi$ and $\Psi$, then 
		$X(\tilde s;M)\equiv0$ follows from taking $n\to\infty$.
	\end{proof}



%
%
%

\bibliographystyle{acm}
\bibliography{bibliography}

%
%
%

\end{document}